\documentclass[12pt]{amsart}
\usepackage{amssymb}
\usepackage{amsmath}
\usepackage{amsrefs}
\usepackage{hyperref}
\usepackage{bbm}
\usepackage{xcolor}

\definecolor{refkey}{gray}{.5}   % graylevel for refs
\definecolor{labelkey}{gray}{.5} % graylevel for labels
\usepackage{graphicx}
\usepackage{tikz-cd}
\usepackage{epsfig}
\usepackage{mathtools,amsfonts}
\usepackage{here}%added by Michi
\usepackage{caption}%added by Michi
\usepackage{subcaption}%added by Michi

\usepackage{comment}
\theoremstyle{plain}

\newtheorem{theorem}{Theorem}[section]
\newtheorem{lemma}[theorem]{Lemma}
\newtheorem{sublemma}{Sublemma}

\newtheorem{proposition}[theorem]{Proposition}

\theoremstyle{definition}
\newtheorem{definition}[theorem]{Definition}
\newtheorem{remark}{Remark}

\newtheorem*{question*}{Question}%added by Michi

\numberwithin{equation}{section}

\renewcommand{\phi}{\varphi}

\usepackage{amsmath}

\DeclareMathOperator{\diam}{diam}

\DeclareMathOperator{\dist}{dist}

\allowdisplaybreaks%added by Michi

\title[Coexistence of divergence and convergence phenomena]{On the coexistence of divergence and convergence phenomena for the Fourier-Haar series for non-negative functions}
\thanks{
The first author is partially supported by Japan Society for the Promotion of Science (JSPS) KAKENHI Grant Number 19K03558. 
The second author is supported by the Knut and Alice Wallenberg foundation of Sweden (KAW)
}

\date{\today}

\begin{document}

\author{Michihiro Hirayama}
\author{Davit Karagulyan}

\maketitle
\begin{abstract}

Let $\{H_{n,m}\}_{n,m\in \mathbb{N}}$ be the two dimensional Haar system and $S_{n,m}f$ be the rectangular partial sums of its Fourier series with respect to some $f\in L^1([0,1)^2)$. Let $\mathcal{N}, \mathcal{M}\subset \mathbb{N}$ be two disjoint subsets of indices. We give a necessary and sufficient condition on the sets $\mathcal{N}, \mathcal{M}$ so that for some $f \in L^1([0,1)^2)$, $f \geq 0$ one has for almost every $z\in [0,1)^2$ that
$$
\lim_{n,m \rightarrow \infty;n,m \in \mathcal{N}}S_{n,m}f(z)=f(z)\quad \text{ and }\quad 
\limsup_{n,m \rightarrow \infty;n,m \in \mathcal{M}}|S_{n,m}f(z)|=\infty.
$$
The proof uses some constructions from the theory of low-discrepancy sequences such as the van der Corput sequence and an associated tiling of the plane. This extends some earlier results.
%study the case when the random variables %$\{\omega_n\}$ are absolutely continuous. 
%We show that if for the sequence  there is no %Dvoretzky covering for the sequence $\zeta_n$, then %there can't be covering with respect to any  %absolute continuous density. 
\end{abstract}

\section{Introduction}

Let $\Psi=\{\psi_k\}_{k \in \mathbb{Z}^d}$ with $\psi_k \in L^2(\mathbb{T}^d)$ be an orthonormal system (i.e., $\|\psi_k\|_2=1$ and $\langle \psi_k,\psi_\ell \rangle=0$, when $k \neq \ell$), and $f \in L^1(\mathbb{T}^d)$. 
We consider the rectangular partial sums of the Fourier series with respect to the system $\Psi$, i.e. for every $n=(n_1, \dots, n_d)\in \mathbb{N}^d$, 
\begin{equation} \label{partial}
    S_{n}f=\sum _{\substack{k=(k_1, \dots, k_d)\in \mathbb{Z}^d;\\ |k_i|\leq n_i}}\langle f,\psi_{k}\rangle \psi_{k}. 
\end{equation}
One can as well use other summation methods, however in this paper we will consider only rectangular summation methods.
It is well known that for certain orthonormal systems there exists $f \in L^1(\mathbb{T})$ so that $S_{n}f$ diverges almost everywhere. For instance the classical example by A. Kolmogorov \cite{Kol1923} shows this for the one dimensional trigonometric system. In \cite{Gosselin}, Gosselin proved that for every increasing sequence of natural numbers $\left(n_{k}\right)$ there exists a function $f \in L^1([0,2 \pi))$ such that
$
\sup _{k \in \mathbb{N}}\left|S_{n_{k}}(f)(x)\right|=\infty.
$
Similar functions can also be constructed for the Fourier-Walsh system.
Another classical system, for which divergence phenomena occurs, is the Haar wavelet which will be defined in detail in Appendix \ref{app}. In this paper we are interested in the following question: for $n=(n_1, \dots, n_d)\in \mathbb{N}^d$ define $|n|=\min_{i}n_i$, then
\begin{question*}
Let { $\mathcal{N},\mathcal{M}\subset \mathbb{N}^d$} be two infinite subsets of indices. 
Under which conditions on the sets $\mathcal{N}, \mathcal{M}$ there exists a function $f\in L^1(\mathbb{T}^d)$, with $f \geq 0$, such that 
\[
\lim_{\substack{|n| \rightarrow \infty \\{n \in \mathcal{N}}}}S_{n}f(z)=f(z)\quad \text{ and }\quad 
\limsup_{\substack{|n| \rightarrow \infty \\n \in \mathcal{M}}}|S_{n}f(z)|=\infty
\]
for Lebesgue almost every $z\in \mathbb{T}^d$?
\end{question*}

\begin{comment}
   
$\chi_n$ of this system are defined on the interval [0,1] as follows:
$$
\chi_{1}(t) \equiv 1 \quad \text { on }[0,1]
$$
if $n=2^{m}+k, k=1, \ldots, 2^{m}, m=0,1, \ldots$, then
$$
\chi_{n}= \begin{cases}\sqrt{2^{m}} & \text { for } t \in\left(\frac{2 k-2}{2^{m+1}}, \frac{2 k-1}{2^{m+1}}\right) \\ -\sqrt{2^{m}} & \text { for } t \in\left(\frac{2 k-1}{2^{m+1}}, \frac{2 k}{2^{m+1}}\right) \\ 0 & \text { for } t \notin\left(\frac{k-1}{2^{m}}, \frac{k}{2^{m}}\right)\end{cases}
$$
\end{comment}

In this paper we give a complete answer to this question for the univariate Haar system. We will consider the case {$d=2$}. The corresponding problem for spherical summation methods and for systems such as the trigonometric and Walsh systems appears to be open.

Denote the univariate Haar system by $\mathcal{H}:=\left\{H_m\right\}_{m\in \mathbb{N}^2}$. (see Appendix \ref{app} for definition and properties). 
For $n,m\in \mathbb{N}$, consider the rectangular partial sums of the Fourier-Haar series as in \eqref{partial}
\[
S_{n,m}f=\sum_{\substack{1\leq i \leq n, \\ 1 \leq j \leq m}}\left\langle f, H_{i,j}\right\rangle H_{i,j}.
\]
It is well known that the correct Orlicz class of convergence for this sums is $L^1\ln^+ L^1$ (see \cite{Jessen-Marcinkiewicz-Zygmund1935}, \cite{Saks1934}). Hence, there exist a function $f \in L^1([0,1)^2)$ for which $S_{n,m}f$ diverges almost everywhere. 

%For $n,m \in \mathbb{N}\setminus \{1\}$, let $n=2^k+i$ and $m=2^s+j$, where $k,s\in \mathbb{N}\cup \{0\}$ and $i=1,\dots, 2^k$, $j=1,\dots, 2^s$. 
For $n\in \mathbb{N}\setminus \{1\}$, one can let $n=2^k+i$, where $k\in \mathbb{N}\cup \{0\}$ and $i=1,\dots, 2^k$. 
Given $\mathcal{N},\mathcal{M}\subset \mathbb{N}$, denote 
\begin{equation} \label{def:C}
C=\left\{k\in \mathbb{N}\cup \{0\} \colon 2^k+i\in \mathcal{N} \text{ for some }i\in \{1,\dots ,2^k\}\right\}
\end{equation}
and
\begin{equation} \label{def:D}
D=\{s\in \mathbb{N}\cup \{0\} \colon 2^s+i\in \mathcal{M} \text{ for some }i\in \{1,\dots ,2^s\}\},
\end{equation}
respectively. 

We have the following theorem:

\begin{theorem}\label{main1}
Let $\mathcal{N},\mathcal{M}\subset \mathbb{N}$ be two subsets of indices and let $C,D$ be defined as above. For the Fourier-Haar series there exists a \emph{non-negative} function $f \in L^1([0,1)^2)$ such that for almost every $z\in [0,1)^2$ we have 
\[
\lim_{\substack{n,m \rightarrow \infty; \\n,m\in \mathcal{N}}}S_{n,m}f(z)=f(z),
\]
and
\[
\limsup_{\substack{n,m \rightarrow \infty; \\n,m\in \mathcal{M}}}|S_{n,m}f(z)|=\infty {\color{blue} }
\]
if and only if
\begin{equation} \label{cond:dst}
\sup_{n \in D}\dist(n, C)=\infty,
\end{equation}
where $\dist(n, C)=\min_{m \in C}|n-m|$. \end{theorem}

\begin{comment}
   
We next consider a theorem that is equivalent to the above theorem.
And also
$$
S_{n,m}(x,y)= \begin{cases}\frac{1}{|I_{i,+}||I_{j,+}|}\int_{I_{i,+}\times I_{j,+}}f(t_1,t_2)dt_1 dt_2 & \text { if } (x,y)\in I_{i,+}\times I_{j,+}\\
\frac{1}{|I_{i,-}||I_{j,+}|}\int_{I_{i,-}\times I_{j,+}}f(t_1,t_2)dt_1 dt_2 & \text { if } (x,y)\in I_{i,-}\times I_{j,+}\\ 
\frac{1}{|I_{i,+}||I_{j,-}|}\int_{I_{i,+}\times I_{j,-}}f(t_1,t_2)dt_1 dt_2 & \text { if } (x,y)\in I_{i,+}\times I_{j,-}\\ 
\frac{1}{|I_{i,-}||I_{j,-}|}\int_{I_{i,-}\times I_{j,-}}f(t_1,t_2)dt_1 dt_2 & \text { if } (x,y)\in I_{i,-}\times I_{j,-}\end{cases}
$$

In this paper we are interested in the following question. Assume we are given two sets $\mathcal{C}, \mathcal{D}\subset \mathbb{N}$, with $\mathcal{D}\cap \mathcal{C}= \emptyset$. 
For $\mathcal{A}\subset \mathbb{Z}$ We also define 
\[
\delta _{\Phi}(z, \mathcal{A},f)=\limsup _{n,m \rightarrow \infty; n,m \in \mathcal{A}}\left|\mathcal{S}_{n,m}f(z)-f(z)\right|.
\]

That is want that Fourier series to be convergent along the sub-sequence $\mathcal{C}$
In this paper we answer this question with for the Fourier-Haar system. 

\end{comment}

We now state an alternative formulation of the above theorem.
Let $\mathcal{R}$ be the family of half-closed axis-parallel rectangles in $\mathbb{R}^{2}$, i.e. $\mathcal{R}=\{[a, b) \times[c, d)\colon a<b, c<d\}$.
For $R \in \mathcal{R}$ we denote by $\diam R$ the length of the diagonal of $R$.
Then let $\mathcal{R}^{\textrm{dyadic}}$ be the family of all dyadic rectangles in $[0,1)^2$ of the form
\[
R_{ n, m}(i,j)=\left[\frac{i-1}{2^{n}}, \frac{i}{2^{n}}\right) \times \left[\frac{j-1}{2^{m}}, \frac{j}{2^{m}}\right), 
\]
where $n, m \in \mathbb{N}$, $1\leq i\leq 2^n$, and $1\leq j \leq 2^m$.

\begin{definition}
A family of rectangles $\mathcal{F} \subset \mathcal{R}$ is said to be a basis of differentiation (or simply a basis), if for any point $z \in \mathbb{R}^{2}$ there exists a sequence of rectangles $R_{k} \in \mathcal{F}$ such that $z \in R_{k}$, $k \in \mathbb{N}$, and $\textrm{diam} R_{k} \rightarrow 0$ as $k \rightarrow \infty $.
\end{definition}

For an infinite subset of integers $C \subset \mathbb{N}$ one can generate a rare basis as follows
\[
\mathcal{F}_{C}=\left\{R_{n, m}(i,j)\colon n,m \in C, i\in \{1,\dots, 2^n\}, j \in \{1,\dots,2^m\} \right\}.
\]

Let $\mathcal{F} \subset \mathcal{R}$ be a differentiation basis. 
For any function $f \in L^{1}\left(\mathbb{R}^{2}\right)$ and $z\in \mathbb{R}^2$ we define
\[
\delta_{\mathcal{F}}(z, f)=\limsup _{\substack{\diam R \rightarrow 0; \\z \in R \in \mathcal{F}}}\left|\frac{1}{|R|} \int_{R} f\, d xdy-f(z)\right|.
\]
(Here and below, let $|\cdot|$ denote the Lebesgue measure on $\mathbb{R}^2$.) 
The function $f \in L^{1}\left(\mathbb{R}^{2}\right)$ is said to be differentiable at a point $z \in \mathbb{R}^{2}$ with respect to the basis $\mathcal{F}$, if $\delta_{\mathcal{F}}(z, f)=0$.

We now state another theorem:

\begin{theorem} \label{thm:main2}
Let $C,D\subset \mathbb{N}$ be two infinite subsets of integers and let $\mathcal{F}_{C}$ and $\mathcal{F}_{D}$ be the corresponding basis. Then there exists a function $f\in L^1([0,1)^2)$, with $f \geq 0$, such that for almost every $z\in [0,1)^2$ we have
\begin{equation*} \label{C1}
    \delta _{\mathcal{F}_C}(z, f)=0,
\end{equation*}
and
\begin{equation*} \label{C2}
    \delta _{\mathcal{F}_D}(z, f)=\infty
\end{equation*}
if and only if
\begin{equation*}\label{dst}
\sup_{n \in D}\dist(n, C)=\infty.
\end{equation*}
\end{theorem}

Let $\mathcal{F}$ be a differentiation basis and consider classes of functions
\begin{align*}
\mathcal{L}(\mathcal{F})&=\left\{f \in L^1\left(\mathbb{R}^{2}\right)\colon \delta_{\mathcal{F}}(z, f)=0 \text { for almost every } z\right\}, \\
\mathcal{L}^{+}(\mathcal{F})&=\left\{f \in L^1\left(\mathbb{R}^{2}\right)\colon \delta_{\mathcal{F}}(z, f)=0 \text{ and } f\geq 0 \text { for almost every }z\right\}.
\end{align*}
Note that $\mathcal{L}(\mathcal{F})$ is the family of functions having almost everywhere differentiable integrals with respect to the basis $\mathcal{F}$.

Note that it is known that \cite{Zerekidze1985}
\[
\mathcal{L}^{+}\left(\mathcal{R}^{\text{dyadic }}\right)=\mathcal{L}^{+}(\mathcal{R}).
\]
This means that for positive functions the basis $\mathcal{R}$ is equivalent to the basis of all dyadic rectangles $\mathcal{R}^{\text{dyadic}}$. 
Remark, however that we do not have $\mathcal{L}\left(\mathcal{R}^{\text{dyadic }}\right)=\mathcal{L}(\mathcal{R})$, i.e. unlike the class of all non-negative functions, there is no equivalence between the differential basis of all rectangles and the class of dyadic rectangles in the sense that convergence with respect to $\mathcal{R}^{\text{dyadic}}$ does not guarantee convergence with respect to $\mathcal{R}$ and the divergence with respect to $\mathcal{R}$ does not guarantee divergence with respect to $\mathcal{R}^{\text{dyadic}}$.

Note that in \cite{Stokolos2006} the authors prove the above theorem for the case $C=\emptyset$ and an arbitrary, infinite subset $D\subset \mathbb{N}$. We remark that the function constructed in the paper is positive. 

In \cite{Karagulyan-Karagulyan-Safaryan2017} the authors consider the case $C=\mathbb{N}\setminus D$, where $D$ is an arbitrary infinite subset and give a necessary and sufficient condition for the existence of a function $f$. The function constructed in the paper is unbounded both from above and below (hence is not positive). 

In \cite{HK2021} the authors study the problem for the basis $\mathcal{R}$, i.e. for the class of all rectangles. They considered two sets $\mathcal{C}, \mathcal{D}\subset [0,1]$ and give conditions on the sets, under which one can construct a functions $f$ which is convergent with respect to rectangles with sides in $\mathcal{C}$ and divergent with respect rectangles with sides in $\mathcal{D}$. The function constructed in the paper is unbounded both from above and below, hence is not positive. Non-positivity of the function is crucial in the proof. As was mentioned above we have $\mathcal{L}\left(\mathcal{R}^{\text{dyadic }}\right)\neq\mathcal{L}(\mathcal{R})$. Due to the non-constructive nature of the argument the convergence and divergence properties of the function on the rectangles from the bases $\mathcal{F}_D$ is not clear. To overcome this issue a new, constructive approach is needed to the problem. We provide such an approach in this paper.

We now sketch the idea of the proof of Theorem \ref{thm:main2}. The idea is to construct an intermediate function satisfying properties in Proposition \ref{prop:main2}. In order to do so we choose rectangles from the basis $\mathcal{F}_D$ and distribute them in a way that they cover a substantial portion of the unit square, then we distribute the support of $f$ in such a way that the integral averages with respect to each rectangle is larger than the prescribed number $M$ thus full-filling condition \eqref{abv}. Hence, for any point that belongs to any of the rectangles the integral averages will be large.

However, at the same time the distribution of the support of $f$ needs to be such that the integral averages with respect to rectangles from $\mathcal{F}_C$ are small (property \eqref{bll}). If we think of the support of $f$ as being concentrated at finite number of points and assume that each point has the same mass, then the question of estimating the expressions $(1/|R|)\int_R f\, dxdy$ will boil down to computing the number of point-supports that fall inside $R$. This is nothing else but a discrepancy estimates for the rectangles in $\mathcal{F}_C$ and for the set $S=\left\{x_{1}, \ldots, x_{N}\right\}$ that carries the support of $f$. 
Namely, if $\mathcal{S}$ is a collection of axis-parallel rectangles, then one defines discrepancy, using Kuipers and Niederreiter's notation [\cite{Kuipers-Niederreiter1974}, page 93], as follows
\begin{equation}\label{discrp}
D_{N}(S)=\sup _{P\in \mathcal{S}}\left|\frac{\#\{x\in S\colon x\in P\}}{N}- |P|\right|.
\end{equation}

It is known that for the van der Corput sequence the above expressions reaches the lowest possible asymptotic bound, i.e.
\[
D_{N}(S)\leq C\frac{\log N}{N}.
\]
Therefore, it is natural to use this sequence to minimize the discrepancy of the distribution of the support of $f$. We remark that the situation is in fact more complicated than the one described above, however the general idea is the same. 
A natural question arises whether the ideas in this paper can be used to construct a sequence for which its discrepancy with respect to one bases of rectangle is different from its discrepancy with respect to another bases, i.e. $\mathcal{F}_C$ and $\mathcal{F}_D$?

To this end, in Section \ref{sec:vdCtiling} we introduce a tilling that follows the dichotomy of the van der Corput sequence and describes a way of distributing rectangles inside a unit square. 
In Section \ref{subsec:pairing} we consider several van der Corput tilings and create pairings between them which eventually leads us to define the function $f$ in Section \ref{subsec:function}. The support of $f$ is placed at the intersection of the rectangles that are paired with each other.
The resulting function turns out to satisfy the desired properties of Proposition \ref{prop:main2}. Due to the constructive nature of the function we are also able to deduce all the necessary information for the rectangles in the bases $\mathcal{F}_C$ and $\mathcal{F}_D$.

To prove the sufficiency, we note that the maximal function \[
M_1f(z) = \sup_{z \in R \in \mathcal{F}_D} \frac{1}{|R|}\int_R f\, dxdy
\]can essentially be estimated from above by the maximal function 
\[
M_2f(z) = \sup_{z \in R \in \mathcal{F}_C} \frac{1}{|R|}\int_R f\, dxdy,
\] 
if the sets $C$ and $D$ are close.

\subsection{Notations}

Throughout the paper, let $\pi _x\colon \mathbb{R}^2\ni (x,y)\mapsto x\in \mathbb{R}_x$ denote the projection onto $x$-axis and respectively, $\pi _y\colon \mathbb{R}^2\ni (x,y)\mapsto y\in \mathbb{R}_y$ denotes the projection onto $y$-axis. 
Let $|\cdot |$ or $\textrm{Leb}$ denote the Lebesgue measure on $\mathbbm{R}^d$ for $d=1,2$. 
For a finite set $A$, let $\# A$ denote the cardinality.

\section{Auxiliary constructions: the Van der Corput sequence and a tiling of the unit square} \label{sec:vdCtiling}

In this section we make some preparational work for proving Proposition \ref{prop:main1} in Section \ref{sec:mainprop}.
%Given $N\in\mathbb{N}$, let $P=\{ p_0,p_1,\dots ,p_{N-1}\}\subset [0,1)^2$ be the van der Corput set. 
%This set is defined as follows. 
For $i \in \mathbb{N}\cup \{0\}$, let $i=a_0+2a_1+2^2a_2+\cdots $, where $a_j\in \{ 0,1\}$, be the binary expression.
Set 
\[
v(i)=\frac{a_0}{2} +\frac{a_1}{2^2} +\frac{a_2}{2^3} +\cdots .
\]
Then define 
\[
p_i=\left( \frac{i}{N} ,v(i)\right) \in [0,1)^2
\]
for $i=0,1,\dots ,N-1$. 
The set $P=\{ p_0,p_1,\dots ,p_{N-1}\}\subset [0,1)^2$ is called the van der Corput set. 
See \cite{Kuipers-Niederreiter1974}. As was already mentioned in the introduction, the van der Corput sequence is known to have low discrepancy. 
\begin{comment}
   
Namely, if we define the discrepancy of a set $S=\left\{x_{1}, \ldots, x_{N}\right\}$ is defined, using Niederreiter's notation, as
\begin{equation}\label{discrp}
D_{N}(S)=\sup _{ u_1,u_2\leq 1
}\left|\frac{\#\{p\in S :p\in [0,u_1)\times[0,u_2)\}}{N}- u_1 u_2\right|.
\end{equation}

It is known that for the van der Corput sequence one has the following bound
$$
D_{N}(P)\leq C\frac{\log N}{N}.
$$
This is the best bound possible.
\end{comment}

We are given a rectangle $R=[0,b)\times [0,a)$ and $(\xi ,\eta )\in [0,1)^2$. 
Henceforth, to simplify the exposition, we will say that \textit{$R$ is placed at $(\xi ,\eta )$} when we translate $R$ by the vector $(\xi ,\eta )$: 
\[
R\mapsto R(\xi ,\eta )=R+(\xi ,\eta )=[\xi ,b
+\xi )\times [\eta ,a+\eta ). 
\]
Note that $(\xi ,\eta )$ specifies the lower left corner of $R(\xi ,\eta )$. 

Let $\mathcal{D} _1=\{ 1/2^k\in [0,1)\colon k\in \mathbb{N} \}$. 
Given $a,b \in \mathcal{D}_1$, let $R$ be an axis-parallel rectangle with height $a$ and width $b$. 
Let $N=1/(ab)$, and $P=\{ p_0,p_1,\dots ,p_{N-1}\}$ be the van der Corput set. 
We will define a tiling on $[0,1)^2$ which is generated by $R$ and associated with $P$. 
%define the van der Corput tiling on $[0,1)^2$ as follows. 
%Suppose $y=1/2^k$. 
%Let
%\[
%\tau _a=\inf \{i\in \mathbb{N} \colon v(i)\in [0,a)\}. 
%\]
%{\color{blue}this is a bit too fast. Why do we have that $\tau _a=a^{-1}$?} Then one can see that $\tau _a=a^{-1}$, and $v(0), v(1),\dots ,v(\tau _a-1)\in [0,1)$ are distributed equidistantly with intervals of $a$. (For $a=1/2^k$, we may write $\tau _a=\tau _k$ by an abuse of notation.)
First, we place $R$ at $p_0, p_1,\dots ,p_{a^{-1}-1}$, respectively.  
Since $v(0),v(1),\dots ,v(a^{-1}-1)$ are distributed equidistantly with intervals of $a$, the rectangles $R(p_0),R(p_1),\dots ,R(p_{a^{-1}-1})$ are disjoint and 
\[
\pi _y\left( \bigsqcup _{i=0}^{a^{-1}-1}R(p_i)\right) =[0,1),
\]
This finishes the first ``column" of tiling. 

%\begin{definition} \label{def:proj}
%Let $\pi _y\colon \mathbb{R}^2\ni (x,y)\mapsto y\in \mathbb{R}_y$ denote the projection onto $y$-axis and respectively, $\pi _x\colon \mathbb{R}^2\ni (x,y)\mapsto x\in \mathbb{R}_x$ denotes the projection onto $x$-axis. 
%\end{definition}

\begin{figure}
\begin{center}
\begin{subfigure}{0.45\textwidth}
\includegraphics[pagebox=artbox,width=\textwidth]{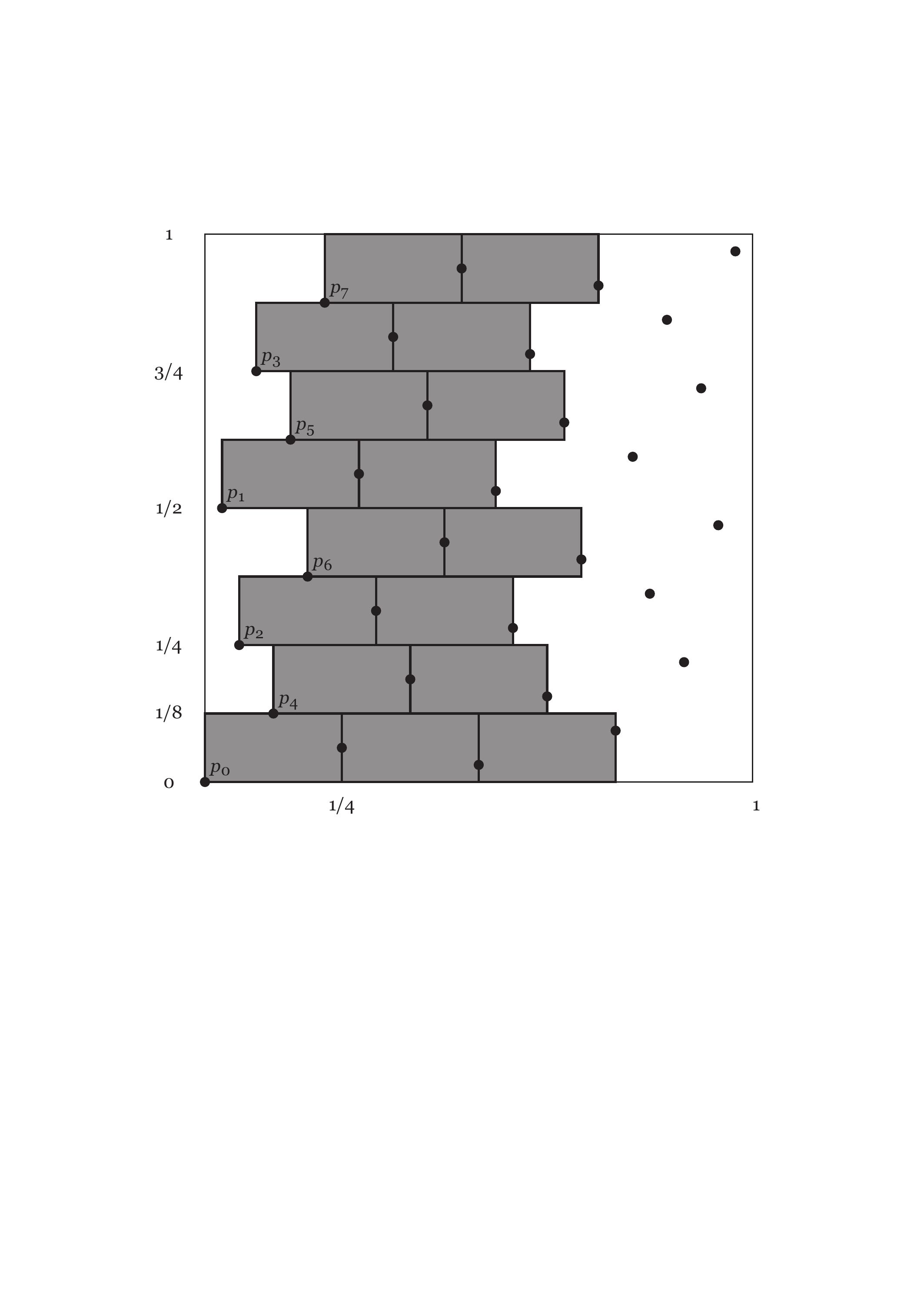}
\caption{The Van der Corput sequence for $N=32$ and the associated Van der Corput tiling for the rectangle with sides $1/8$ and $1/4$.}
\label{fig:vdc_tiling32h}
\end{subfigure}
\hfill
\begin{subfigure}{0.45\textwidth}
\includegraphics[pagebox=artbox,width=\textwidth]{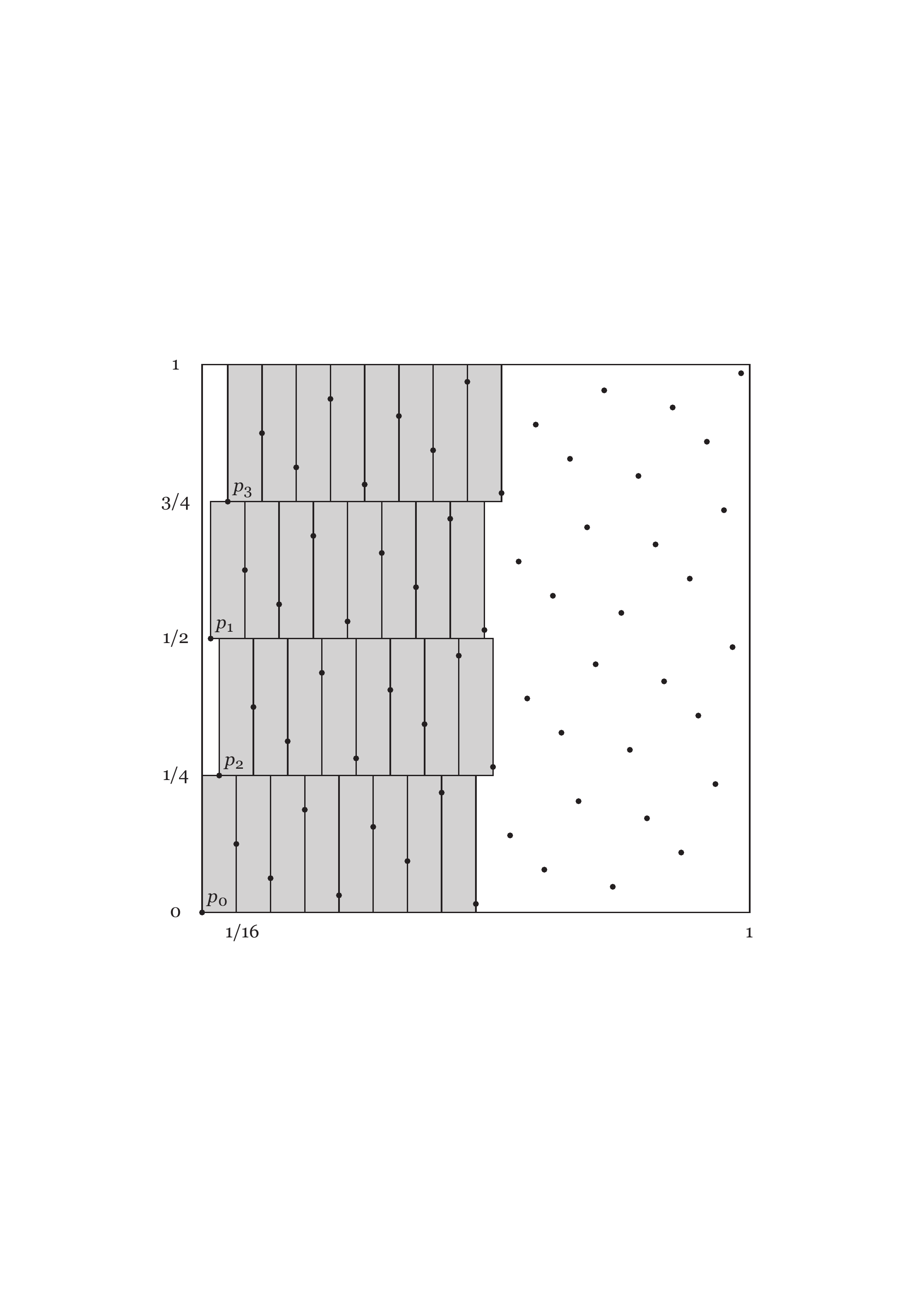}
\caption{The Van der Corput sequence for $N=64$ and the associated Van der Corput tiling for the rectangle with sides $1/4$ and $1/16$.}
\label{fig:vdc_tiling64}
\end{subfigure}
\caption{The $N$-points Van der Corput sequence and the associated Van der Corput tiling. Remark that each tile contains only one van der Corput point.}
\label{fig:vdctiling}
\end{center}
\end{figure}

To determine the second column we now translate the first column by the horizontal vector $(b,0)$ and subsequently by vectors $j(b,0)$, $j=1, \dots, \frac{1}{b}-1$. 
That is we will have the following collections $\{R(p_i)+(jb,0)\colon i\in \{ 0,1,\dots ,a^{-1}-1\}\}$ for each $j\in {0,1,\dots ,b^{-1}-1}$. 
One then can see that the resulting placement of figures will look like Figure \ref{fig:vdctiling}. %Note that in Figure 1 $\tau_a=4$
\begin{comment}
   To determine the next column, we place $R$ at 
\[
\left( \pi _x(p_{\tau _a}),v(0)\right), \left( \pi _x(p_{\tau _a+1}),v(1)\right) , \dots ,\left( \pi _x(p_{2\tau _a-1}),v(\tau _a-1)\right) ,
\]
respectively, 
Note here that 
\[
\pi _x(p_{\tau _a+i})=\frac{\tau _a+i}{N} =\frac{1}{aN} +\frac{i}{N} =b+\frac{i}{N} ,
\]
and hence one sees that 
\[
R(\pi _x(p_{\tau _a+i}),v(i))=R(p_{i})+(b,0)
\]
for $i=0,1,\dots ,\tau _a-1$.
In other words, one sets $\tau _a$-many rectangles $R$ in the same way as in the first but translation by $(b,0)$ along the horizontal direction. 
We repeat $1/b$ times the procedure above. 
\end{comment}
Identifying $\{ 0\} \times [0,1)\ni (0,a)\sim (1,a)\in \{ 1\} \times [0,1)$ will give a tiling of $[0,1)^2$ generated by $R$. 
We denote the collection of all tiles by $\mathcal{T}_{a,b}$, more specifically
\[
\mathcal{T}_{a,b}=\{ R(p_{i})+(jb,0)\colon i\in \{ 0,1,\dots ,a^{-1}-1\} ,j\in {0,1,\dots ,b^{-1}-1}\} ,
\]
and thus $\# \mathcal{T}_{a,b}=a^{-1}\times b^{-1}=N=\# P$. 

One can see from the figures in Figure \ref{fig:vdctiling} that each horizontal row of rectangles is a horizontal translation of other rows. 
Therefore, any row can be described by the amount of horizontal translation vector with respect to the bottom row $[0,b)\times [0,a) + j(b,0)$, where $j=0,\dots, b^{-1}-1$. 
The horizontal translation length will be denoted by $d_*(\ell)$ with $\ell =0,\dots, a^{-1}-1$, starting from the bottom row. Thus the $\ell$'th row can be given as
\begin{equation}\label{cond:horizontal_transl}
    [0,b)\times [0,a) + j(b,0)+(d_*(\ell),\ell a).
\end{equation}

One can see that the sequence $d_*(\ell)$ is similar to the $y$ coordinate of the van der Corput sequence in the sense that $d_*(\ell)=v(\ell)b$.

\subsection{Pairing} \label{subsec:pairing}

In this Section we consider two collections of van der Corput tilings $\mathcal{T}_{a,b}$ and $\mathcal{T}_{x,y}$, with $x>a$,  $b>y$ and define a pairing between the rectangles from each collection. 

For a rectangle $R=[x_1,x_2)\times [y_1,y_2)$, let $\partial _{\mathrm{left}}R=\{ x_1\} \times [y_1,y_2)$ and $\partial _{\mathrm{v}}R=(\{ x_1\} \times [y_1,y_2))\sqcup (\{ x_2\} \times [y_1,y_2))$. 
%We will say $R$ is well-placed with respect to $t\in \mathcal{T}_{a,b}$ if $\partial _{\mathrm{v}}R\cap \partial _{\mathrm{v}}t\neq \emptyset$. 
%We will say $R$ is well-placed with respect to $\mathcal{T}_{a,b}$ if there is $t\in \mathcal{T}_{a,b}$ such that $\partial _{\mathrm{v}}R\cap \partial _{\mathrm{v}}t\neq \emptyset$. 

\begin{lemma}[Pairing lemma] \label{lem:pairing}
Let $a,b,x,y\in \mathcal{D} _1$ so that $x>a$, $y<b$ and $xy\leq ab$, and $n>5$. 
Then for any $t\in \mathcal{T}_{a,b}$ one can find a collection $\mathcal{P}_{x,y}(t)\subset \mathcal{T}_{x,y}$ consisting of adjacent rectangles so that the following properties hold:
\begin{enumerate}
    \item[1)] Every $u\in \mathcal{P}_{x,y}(t)$ intersects with $t$ so that $\pi_y(u)\supset \pi _y(t)$ and $\pi_x(u)\subset \pi_x(t)$. 
    
    \item[2)] The union $\cup_{u\in \mathcal{P}_{x,y}(t)}u$ is a rectangle of height $x$ and width $ab/x$ such that $\pi_x(t)\setminus \pi_x\left(\cup_{u\in \mathcal{P}_{x,y}(t)}u\right)$ consists of two components (intervals) of length at least $2b/n$, and $\# \mathcal{P}_{x,y}(t)=ab/(xy)$. 
    
    \item[3)] For different $t\in \mathcal{T}_{a,b}$ the corresponding unions (rectangles) $\cup_{u\in \mathcal{P}_{x,y}(t)}u$ have disjoint interiors. However, they may have common boundaries. 
    
    \item[4)] For every $t\in \mathcal{T}_{a,b}$, one has 
    
    \[ \left|\bigcup_{u\in \mathcal{P}_{x,y}(t)}u \right|=|t|. \]
    %\[ 1-\frac{xy}{ab} \leq \frac{|\cup \mathcal{P}_{x,y}^*(t)|}{|t|}\leq 1. \]
    
    \item[5)] \[ [0,1)^2=\bigcup_{t\in \mathcal{T}_{a,b}}  \left(\bigcup_{u\in \mathcal{P}_{x,y}(t)}u\right).\]
    %\[ \left| [0,1)^2\setminus \sqcup_{t\in \mathcal{T}_{a,b}}  (\sqcup \mathcal{P}_{x,y}^*(t))\right| \leq \frac{xy}{ab}.\]
\end{enumerate}
\end{lemma}

\begin{figure}
\begin{center}
\includegraphics[pagebox=artbox,width=12cm]{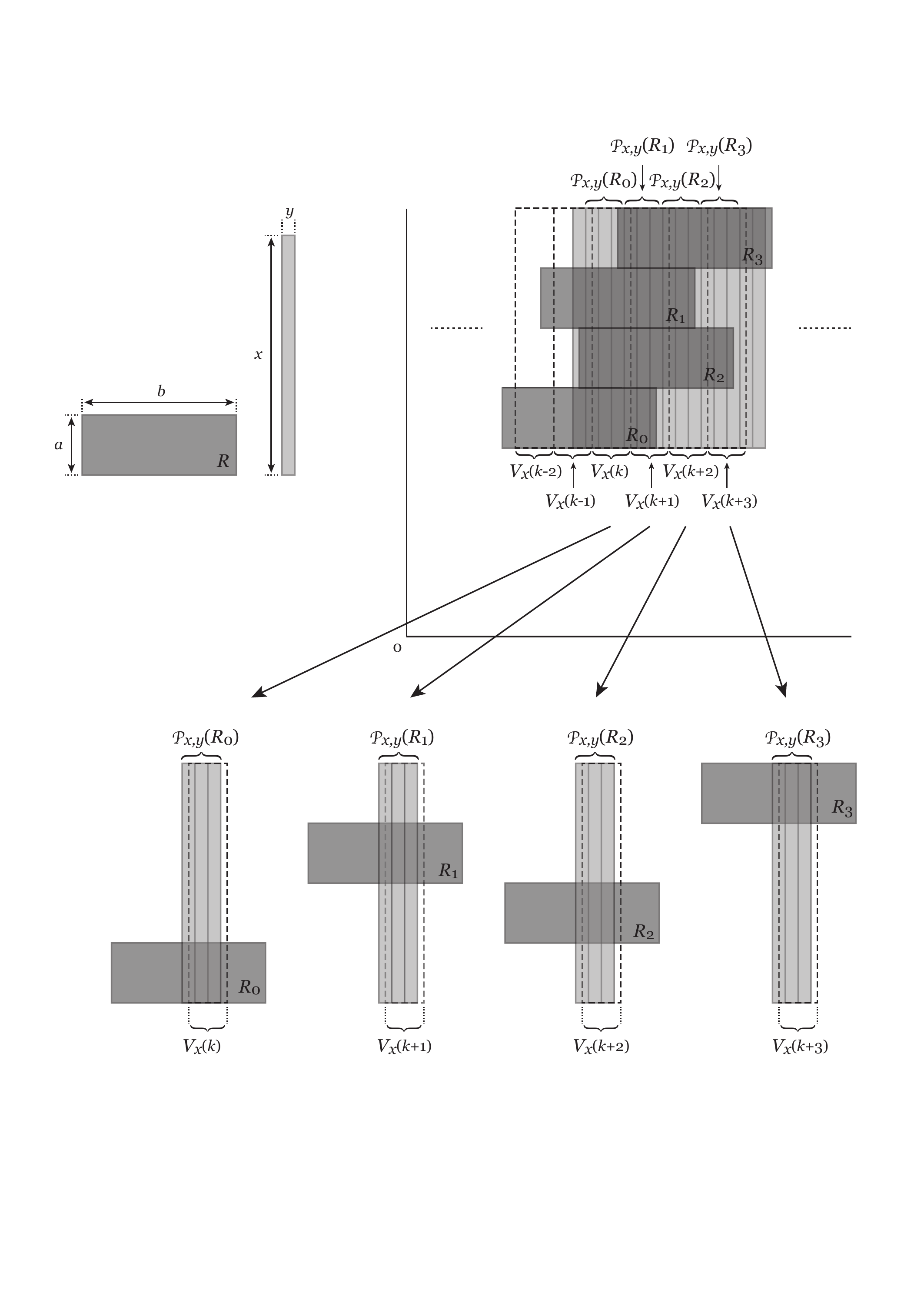}
\end{center}
\caption{Pairing construction for $R\in \mathcal{T}_{a,b}$. The light gray rectangles belong to $\mathcal{T}_{x,y}$. The dashed ones are $V_x$. $k=k_\textrm{m}$ is defined in \eqref{cond:translate_mid}.}
\label{fig:pairing}
\end{figure}

\begin{proof}

Let  
\[
H_x=[0,1)\times [0,x)
=\pi_y^{-1}([0,x))
\]
be a horizontal strip. 
First, we will define $\mathcal{P}_{x,y}(t)$ for $t\in \mathcal{T}_{a,b}$ with $t\subset H_x$.  
For $k=0,1,\dots ,(x/(ab))-1$, let 
\[
V_x(k)=\left[ k\frac{ab}{x}, (k+1)\frac{ab}{x}\right) \times [0,x) \subset H_x.
\]
%Clearly $\sqcup _{j=0}^{(x/a)-1}V_x(j)=[0,b)\times [0,x)$, and 
Note that $|V_x(k)|=ab=|t|$.

Consider $\mathcal{T}_{a,b}(x)=\left\{ R(p_i)\in \mathcal{T}_{a,b}\colon R(p_i)\subset H_x,i\in \{0,1,\dots ,a^{-1}-1\}\right\}$. 
Let $0=i_0<i_1<\dots <i_{(x/a)-1}\leq a^{-1}-1$ be the indices so that $R(p_{i_j}) \in \mathcal{T}_{a,b}(x)$.
(Note that $\mathcal{T}_{a,b}(x)$ is the collection of tiles $R_0, R_1, R_2, R_3$ in Figure \ref{fig:pairing} corresponding to the case when the lower left corner of $R_0$ is the origin with an abbreviation $R_j=R(p_{i_j})$.) 
Next, sort the horizontal drifts for $R(p_{i_j})\in \mathcal{T}_{a,b}(x)$ in an ascending order. 
To do this, let $\sigma \colon \{ 0,1,\dots ,(x/a)-1\} \to \{0,1,\dots, (x/a)-1\}$ be a permutation such that 
\[
R(p_{i_j})\subset [0,1)\times [\sigma (j)a,(\sigma (j)+1)a)
\]
for every $R(p_{i_j})\in \mathcal{T}_{a,b}(x)$. 
Then one has 
\[
0=d_*(\sigma(0))<d_*(\sigma(1))<\dots <d_*(\sigma((x/a)-1)),
\]
where $d_*$ is the horizontal translation length for $\mathcal{T}_{a,b}$. 
%Namely, consider the set $d_*(0),d_*(1),\dots ,d_*((x/a)-1)$ and find a permutation $\sigma \colon \{ 0,1,\dots ,(x/a)-1\} \to \{0,1,\dots, (x/a)-1\}$ such that 
%\[
%0=d_*(\sigma(0))<d_*(\sigma(1))<\dots <d_*(\sigma((x/a)-1)). 
%\]
Note here that
\[
d_*(\sigma (j))%=v(\ell_j)b
=\frac{j/x}{N} =\frac{jab}{x} %=d(p_{j/x})
\]
for $j\in \{ 0,1,\dots ,(x/a)-1\}$. 
It follows that 
\begin{equation} \label{difference_drifts}
    d_*(\sigma (j+1))-d_*(\sigma (j))=\frac{(j+1)ab}{x} -\frac{jab}{x}=\frac{ab}{x}
\end{equation}
for $j\in \{ 0,1,\dots ,(x/a)-1\}$. 
By \eqref{difference_drifts}, for every $j\in \{ 0,1,\dots ,(x/a)-1\}$, one sees that $V_x(j)$ is well-placed with respect to $R(p_{i_j})\in \mathcal{T}_{a,b}(x)$ in the sense that $\partial _{\mathrm{left}} R(p_{i_j})\subset \partial _{\mathrm{left}} V_x(j)$. 
%Here and below, let $\sigma $ denotes the permutation of $\{ 0,1,\dots ,(x/a)-1\}$ such that 
%\[
%R\left( p_{\sigma (\ell _j)}\right) \subset [0,1)\times [ja,(j+1)a)
%\]
%for $j\in \{ 0,1,\dots ,(x/a)-1\}$. 
%Then we associate each $V_x(j)$ with $R(p_{\sigma (\ell _j)})\in \mathcal{T}_{a,b}$. 
%Then we associate $R(p_{\sigma (\ell _j)})$ with the tiles from  $\mathcal{T}_{x,y}$ that fully fall inside $V_x(j)$, that is, 
Let $k\in \mathbb{N}$ such that 
\begin{equation*} \label{cond:translate1}
    \frac{2b}{n} <k\frac{ab}{x} <b-\frac{2b}{n},
\end{equation*}
or equivalently
\begin{equation*} \label{cond:translate2}
    \frac{2x}{an} <k<\frac{(n-2)x}{an}. 
\end{equation*}
For instance, one can choose an integer $k$ so that
\begin{equation}  \label{cond:translate_mid}
    k=\frac{x}{2a}
\end{equation}
for $n>4$. 
To distinguish we denote the $k$ above by $k_\textrm{m}$. 
Then we pair $R(p_{i_j})$ with the tiles from  $\mathcal{T}_{x,y}$ that fully fall inside $V_x(j+k_{\text{m}})$ or intersect $\partial _{\mathrm{left}} V_x(j+k_{\text{m}})$, that is define 
\begin{equation} \label{paired_family}
    \mathcal{P}_{x,y}\left( R\left(p_{i_j}\right)\right) =\{ u\in \mathcal{T}_{x,y} \colon u\subset V_x(j+k_{\text{m}})\} \cup \{ u\in \mathcal{T}_{x,y}\colon u\cap \partial _{\mathrm{left}} V_x(j+k_{\text{m}})\neq \emptyset\}.
\end{equation}
Thus for each $j\in \{ 0,1,\dots ,(x/a)-1\}$, one has $u\cap R(p_{i_j})\neq \emptyset $ for every $u\in \mathcal{P}_{x,y}\left( R(p_{i_j})\right)$. 
Note also that $\mathcal{P}_{x,y}\left( R\left(p_{i_j}\right)\right)$ will be disjoint for distinct $j\in \{ 0,1,\dots ,(x/a)-1\}$. 
Next note that for $R\in \mathcal{T}_{a,b}$, with $R \subset H_x$, we have $\mathcal{P}_{x,y}(R)=\{ u\in \mathcal{T}_{x,y} \colon u\subset V_x(j+k_{\text{m}})\}=V_x(j+k_{\text{m}})$. 

However, for $R\in \mathcal{T}_{a,b}$ with $R\subset [0,1)\times [\ell x,(\ell +1)x)$, $\ell\geq 1$, a tile $u\in \mathcal{T}_{x,y}$ that intersecting $V_x(j+k_{\text{m}})$ partially will exists. 
Hence there will be a miss-match between $\mathcal{P}_{x,y}\left( R\right)$ and $V_x(j+k_m)$ as in Figure \ref{fig:pairing}. %(See Remark \ref{rem:pairing} below.) 
This justifies the definition of $\mathcal{P}_{x,y}\left( R\left(p_{i_j}\right)\right)$ in \eqref{paired_family}. 

One can see that for every $j\in \{ 0,1,\dots ,(x/a)-1\}$, we have $\# \mathcal{P}_{x,y}\left( R\left(p_{i_j}\right)\right) = ab/(xy)$. 
%There will be at most two rectangles from $\mathcal{T}_{x,y}$ that will intersect $V_x(j)$ partially, namely the one's at the vertical boundary  $\partial _{\mathrm{v}}V_x(j)$ which will have total measure at most $2x y$. 
Thus, one has 
%\[ 1-\frac{x y}{\left| R(p_{\sigma (\ell _j))}\right| } \leq \frac{\left| \bigsqcup \mathcal{P}_{x,y}^*\left( R(p_{\sigma (\ell _j)})\right) \right| }{\left| R(p_{\sigma (\ell _j)})\right| } \leq 1 \]
\[
\left| \bigcup_{u\in \mathcal{P}_{x,y}\left( R\left(p_{i_j}\right)\right)}u \right| =\left| R\left(p_{i_j}\right) \right|
\]
for every $j\in \{ 0,1,\dots ,(x/a)-1\}$. 

We now define pairing for remaining tiles from $\mathcal{P}_{a,b}$ that belong to $H_x$. 
We recall that the tiling $\mathcal{T}_{a,b}$ in $H_x$ is just a translations of the rectangles in $\mathcal{T}_{a,b}(x)$ by $(b,0)$. 
Note also that 
\[
V_x(j)+(b,0)=\left[ b+j\frac{ab}{x}, b+(j+1)\frac{ab}{x}\right) \times [0,x)=V_x((x/a)+j)\subset H_x
\]
for $j=0,1,\dots ,(x/a)-1$. 
Hence repeating the argument above, we can associate to every $R\in \mathcal{T}_{a,b}$ in $H_x$ a tile from $\mathcal{T}_{x,y}$ that fully fall inside $V_x(j)$ for some $j\in \{ 0,1,\dots, (x/ab)-1\} $. 
We have defined $\mathcal{P}_{x,y}(t)$ for every $t\in \mathcal{T}_{a,b}$ with $t\subset H_x$. 

Recall that, each horizontal row of tiling  $\mathcal{T}_{a,b}$ in $H_x^{(\ell)}=[0,1)\times [\ell x,(\ell +1)x)$, $\ell =1,\dots, x^{-1}-1$, is a horizontal translation of the tilling in $H_x$. 
To define $\mathcal{P}_{x,y}(t)$ for $t\in \mathcal{T}_{a,b}$ with $t\subset H_x^{(\ell)}$, consider 
\[
V_x^{(\ell)}(k)=\left[ k\frac{ab}{x}, (k+1)\frac{ab}{x}\right) \times [\ell x,(\ell +1)x) +(d_*(\ell x),0)\subset H_x^{(\ell)}
\]
for $k=0,1,\dots ,(x/(ab))-1$, and 
\[
\mathcal{T}_{a,b}^{(\ell)}(x)=\left\{ R(p_i)\in \mathcal{T}_{a,b}\colon R(p_i)\subset H_x^{(\ell)},i\in \{0,1,\dots ,a^{-1}-1\}\right\}.
\]
Here $d_*$ is the horizontal translation length for $\mathcal{T}_{a,b}$. 
In view of \eqref{cond:horizontal_transl}, note that each $R\in \mathcal{T}_{a,b}^{(\ell)}(x)$ can be given as 
\[
R=R_0+\ell (0,x)+(d_*(\ell x),0)
\]
for some $R_0\in \mathcal{T}_{a,b}(x)$. 
Hence, by the same argument for $t\in \mathcal{T}_{a,b}(x)$, one can associate to each $t\in \mathcal{T}_{a,b}^{(\ell)}(x)$ a collection $\mathcal{P}_{x,y}(t)$, the set defined as in \eqref{paired_family}. 
Therefore, repeating the same construction above, one can associate to each $t\in \mathcal{T}_{a,b}$ a collection $\mathcal{P}_{x,y}(t)$, and the properties 1), 2), 3), and 4) will follow. 
Since every $u\in \mathcal{T}_{x,y}$ belongs to a (unique) family $\mathcal{P}_{x,y}(t)$, we have property 5). 
%We consider the pairs $(R,V)$, with $R,V$ from $H_x$. 
%Remark that for every $t\in \mathcal{T}_{a,b}$ we have uniquely defined set $V$ associated to it.
%Therefore, one can associate to each $t\in \mathcal{T}_{a,b}$ a collection $\mathcal{P}_{x,y}(t)$, the set defined as in \eqref{paired_family}, and the properties 1), 2), 3), and 4) follow. 
%Since every $u\in \mathcal{T}_{x,y}$ belongs to a (unique) family $\mathcal{P}_{x,y}(t)$, we have property 5).
%Then it follows from 2) and 3) that \[ \left| [0,1)^2\setminus \sqcup _{t\in \mathcal{T}_{a,b}}\sqcup \mathcal{P}_{x,y}^*(t)\right| \leq \frac{x y}{ab}, \] and which yields property 4). 
\end{proof}

\begin{comment}
   \begin{remark} \label{rem:pairing}
As it was mentioned in the proof above, there will exists a tile $u\in \mathcal{T}_{x,y}$ such that $u\cap \partial _{\mathrm{left}} V_x(j+i_{\text{m}})\neq \emptyset$ for $R\in \mathcal{T}_{a,b}$ in general. 
See \eqref{paired_family}. 
This is due to the fact that the difference of drifts $d_{a,b}-d_{x,y}$ needs not be an integer multiple of $y$ in general. 
See Figure \ref{fig:incomplete_tile}. 
\begin{figure}[H]
\begin{center}
\includegraphics[pagebox=artbox,width=5cm]{pairing1.pdf}
\end{center}
\caption{The shaded rectangle is an ``incomplete" tile $u\in \mathcal{P}_{x,y}^*(R)$ for some $R\in \mathcal{T}_{a,b}$.}
\label{fig:incomplete_tile}
\end{figure}
\end{remark} 
\end{comment}

\subsection{A geometric construction}

Let $b_1>b_2>\dots >b_{2n}$ with $b_i\in \mathcal{D}_1$ and $n>5$ be a power of two. 
Suppose that the sequence $\{ b_k\}$ satisfies 
\begin{equation} \label{cond:decay1}
    b_{n}b_{n+1}>b_{n-1}b_{n+2}>\dots >b_{2}b_{2n-1}>b_{1}b_{2n}. 
\end{equation}
We now construct chains of ``admissible'' tiles from the sequence $\{\mathcal{T}_{b_{n-k}, n b_{n+k+1}}\} _{k=0}^{n-1}$. %, and in the next section we will define those ``admissible'' ones. 
First, consider a truncation of tiles from each $\mathcal{T}_{b_{n-k}, n b_{n+k+1}}$ as follows. 
Given a rectangle $R=[x_1,x_2)\times [y_1,y_2)$, $n\in \mathbb{N}$, and $\omega \in [0,x_2-x_1-\frac{x_2-x_1}{n})$, let 
\begin{equation}\label{def:r}
R^*(\omega) =\left[ x_1+\omega,x_1+\omega+\frac{x_2-x_1}{n}\right) \times [y_1,y_2).
\end{equation}
\begin{comment}
   \begin{sublemma}\label{sb-lem}
   Let $a,b,x,y\in \mathcal{D} _1$ so that $x>a>b>y$, and $n\in \mathbb{N}$ a power of two. 
   Then for $R\in \mathcal{T}_{a,b}$, one has
   \[
   \left| \bigsqcup _{\substack{A\in \mathcal{T}_{x,y}, \\ A\cap R\neq \emptyset }}R'\cap A'\right| \leq \frac{ab}{n^2}+\frac{2a y}{n}.
   \]
   \end{sublemma}
   \begin{proof}
   Since $\# \{ A\in \mathcal{T}_{x,y} \colon A\cap R\neq \emptyset \} \leq \frac{b/n}{y} +2 =\frac{b}{n y} +2$ ({\color{blue} or may be $\leq b/(n y)+1$}), it follows that 
   \[
   \left| \bigsqcup _{\substack{A\in \mathcal{T}_{x,y}, \\ A\cap R\neq \emptyset }}R'\cap A'\right| \leq a\frac{y}{n}\times \left( \frac{b}{n y} +2\right) =\frac{ab}{n^2} +\frac{2a y}{n}.
   \]
   \end{proof}
\end{comment}
%Then to each tiling $\mathcal{T}_{b_{n-k}, n b_{n+k+1}}$ we associate a collection $\mathcal{T'}_{b_{n-k}, n b_{n+k+1}}$ defined as follows
%\[
%\mathcal{Q}_{b_{n-k}, n b_{n+k+1}}=\left\{ \widetilde{R}\colon R\in \mathcal{T}_{b_{n-k}, n b_{n+k+1}} \right\}.
%\]
Note that $R^*(\omega) \subset R$. 
By Lemma \ref{lem:pairing}-5), for every $k\in \{0,1,\dots ,n-1\}$, one has  
\[
\textrm{Leb} \left(\bigcup_{R \in\mathcal{T}_{b_{n-k}, n b_{n+k+1}}}R^*(\omega_{R}) \right) =\frac{1}{n}
\]
for every $\omega_R\in [0,(n-1)b_{n+k+1})$. 
(Note here that $\omega=\omega_R$ can vary with each $R\in\mathcal{T}_{b_{n-k}, n b_{n+k+1}}$.) 
Let $\lambda \in (0,1)$ be such that for all $k \in \{n+1,\dots ,2n-1\}$
\begin{equation} \label{cond:decay2}
    \frac{b_{k+1}}{b_{k}}<\lambda.
\end{equation}

\begin{lemma}\label{Lmm-cov}
Let $\{ b_k\} $ be a sequence as above. 
Suppose \eqref{cond:decay2} holds with $\lambda <\frac{1}{2n(n-1)}$. 
Then for arbitrary choices of translation vectors {$\{\omega_{R}\}$}
\[
\left|\bigcup_{k=0}^{n-1}\bigcup_{R \in\mathcal{T}_{b_{n-k}, n b_{n+k+1}}}
R^*(\omega_{R}) \right|\geq \frac{1}{2}.
\]
\end{lemma}
\begin{proof}
For each $k\in \{ 0,1,\dots ,n-1\}$, let 
\[
W_k=\bigcup_{R \in\mathcal{T}_{b_{n-k}, n b_{n+k+1}}}R^*(\omega_{R}).
\]
Then one has $\textrm{Leb} (W_k)=1/n$, and 
\begin{align*}
    \textrm{Leb} \left( [0,1)^2\setminus \bigcup_{k=0}^{n-1}W_k\right) 
    &=\int \prod_{k=0}^{n-1} (1-\mathbbm{1}_{W_k}(x,y))\, dxdy \\
    &\leq 1 -\sum_{k=0}^{n-1} \int \mathbbm{1}_{W_k}\, dxdy + \sum_{0\leq k < \ell < n }\int \mathbbm{1}_{W_k} \mathbbm{1}_{W_\ell}\, dxdy \\
    &=\sum_{0\leq k < \ell < n}\int \mathbbm{1}_{W_k} \mathbbm{1}_{W_\ell}\, dxdy.
\end{align*}

For $k<\ell$, observe that for a given $R\in \mathcal{T}_{b_{n-k}, n b_{n+k+1}}$, we have
\[
\# \left\{ A\in \mathcal{T}_{b_{n-\ell }, n b_{n+\ell +1}} \colon A\cap R^*(\omega_{R})\neq \emptyset \right\} \leq \frac{n b_{n+k+1}/n}{n b_{n+\ell +1}} +2 =\frac{b_{n+k+1}}{n b_{n+\ell +1}} +2
\]
and 
\[
\textrm{Leb} \left( R^*(\omega_k) \cap A^*(\omega_{A})\right) \leq b_{n-k}\frac{n b_{n+\ell +1}}{n} =b_{n-k}b_{n+\ell +1}
\]
Hence
\begin{align*}
    \int \mathbbm{1}_{W_k} \mathbbm{1}_{W_\ell}\, dxdy 
    &\leq \sum _{R\in \mathcal{T}_{b_{n-k}, n b_{n+k+1}}}\textrm{Leb} \left( \bigsqcup _{\substack{A\in \mathcal{T}_{b_{n-\ell }, n b_{n+\ell +1}}\colon \\ A\cap R^*(\omega_{R})\neq \emptyset }}R^*(\omega_{R})\cap A^*(\omega_{A})\right) \\
    &\leq \sum _{R\in \mathcal{T}_{b_{n-k}, n b_{n+k+1}}} \left( \frac{b_{n+k+1}}{n b_{n+\ell +1}} +2\right) \times  b_{n-k}b_{n+\ell +1} \\
    &=\sum _{R\in \mathcal{T}_{b_{n-k}, n b_{n+k+1}}} \left( \frac{b_{n-k}b_{n+k+1}}{n} +2b_{n-k}b_{n+\ell +1} \right) \\
    &\leq \frac{1}{n b_{n-k}b_{n+k+1}} \left( \frac{b_{n-k}b_{n+k+1}}{n} +2b_{n-k}b_{n+\ell +1} \right) = \frac{1}{n^2} +\frac{2b_{n+\ell +1}}{n b_{n+k+1}} 
\end{align*}
as $\# \mathcal{T}_{b_{n-k}, n b_{n+k+1}}=1/(n b_{n-k}b_{n+k+1})$. 
In view of \eqref{cond:decay2}, it follows that 
\begin{align*}
    \textrm{Leb} \left( [0,1)^2\setminus \bigcup_{k=0}^{n-1}W_k\right) 
    &\leq \sum_{0\leq k < \ell < n}\int \mathbbm{1}_{W_k} \mathbbm{1}_{W_\ell}\, dxdy \\ 
    &\leq \frac{n(n-1)}{2} \left( \frac{1}{n^2} +\frac{2b_{n+\ell +1}}{n b_{n+k+1}} \right) \\
    &\leq \frac{n(n-1)}{2} \left( \frac{1}{n^2} +\frac{2}{n} \lambda \right) <\frac{1}{2} 
\end{align*}
since $\lambda <\frac{1}{2n(n-1)}$.
This implies the lemma. 
\begin{comment}
   Denote by $A_n$ the above union. 
   It is a union of $n$ collection of tilings $A_n=\cup_k B_k$. Let $\chi_k(x,y)=\mathbbm{1}_{B_k}(x,y)$. 
   Then
   \[
   \int\prod_k (1-\chi_{k}(x,y))\, dxdy\leq 1 - \int\sum_k \chi_k + \sum_{i \neq j}\int \chi_i \chi_j=\sum_{i \neq j}\int \chi_i \chi_j \, dxdy.
   \]
   Note that for $i \neq j$ due to Sublemma \ref{sb-lem} we have
   $$
   \int \chi_i(x,y) \chi_j(x,y) \, dx dy\lesssim \frac{1}{n^2}.
   $$
   Hence
   $$
   \int\prod_k (1-\chi_{I}(x,y))\, dxdy\leq \frac{n(n-1)}{2}\frac{1}{n^2}\leq \frac{1}{2}.
   $$
   This implies the Lemma.
\end{comment}
\end{proof}

In what follows, we suppose that $\{ b_k\} _{k=1}^{2n}$ is a decreasing sequence satisfying \eqref{cond:decay1} and \eqref{cond:decay2} with $\lambda <1/2n(n-1)$. 
We assume further that 
\begin{equation} \label{cond:decay3}
    b_k<\lambda b_{k-1}\left( <\frac{b_{k-1}}{2n}\right)
\end{equation}
for every $k\in \{ 2,3,\dots ,n\}$. 
For simplicity of notation, we may sometimes write $\mathcal{T}_{b_{n-k}, n b_{n+k+1}}$ by $\mathcal{T}_{k}$ for $k=0,1,\dots ,n-1$. 
%To each family $\mathcal{T}_{k}$ we associate a collection $\mathcal{Q}_{k}$ defined as follows \[ \mathcal{T}_{k}'=\left\{ R'\colon R\in \mathcal{T}_{k} \right\}. \]
Using Lemma \ref{lem:pairing}, we now construct chains of tiles from $\{ \mathcal{T}_{k}\}_{k=0}^{n-1}$. 
Given a tile $t_{0}\in \mathcal{T}_{0}$, using Lemma \ref{lem:pairing} recursively, one can associate to it a collection $\mathcal{P}_{1}(t_{0})\subset \mathcal{T}_{1}$, and for every tile $t_{1}\in \mathcal{P}_{1}(t_{0})$, one can associate $\mathcal{P}_{2}(t_{1})\subset \mathcal{T}_{2}$, and so forth.

\begin{lemma} \label{lem:translate0}
One can choose $\omega =\omega_{t_{0}}\in [0,(n-1)b_{n+1})$ so that every $t_{1}\in \mathcal{P}_{1}\left( t_{0}\right)$ intersects with $t_{0}^*(\omega )$ (see \eqref{def:r} for definition of $t_{0}^*$) in such a way that 
\begin{equation*} %\label{cond:translate0}
   \begin{cases}
    \pi_x\left( t_{1}\right) \subset \pi_x\left( t_{0}^*(\omega )\right), \\
    \pi_y\left( t_{1}\right) \supset \pi_y\left( t_{0}^*(\omega )\right). 
   \end{cases}
\end{equation*}
\end{lemma}
\begin{proof}
Write $[x,x+n b_{n+1})\times [y,y+b_n)$ for $t_{0}\in \mathcal{T}_{0}$. 
Taking $\omega =(n-1)b_{n+1}/2$ implies 
\[
x+\omega +\frac{b_{n+1}}{2} =x+\frac{n b_{n+1}}{2},
\]
and hence $t_{0}^*(\omega )\subset t_{0}$ is placed at the middle of $t_{0}$.  
Since 
\[
x+\frac{n b_{n+1}}{2} \in \pi_x\left( \bigcup_{u\in  \mathcal{P}_{1}\left( t_{0}\right)} u\right) 
\]
by Lemma \ref{lem:pairing}-2), with the aid of \eqref{cond:translate_mid}, it follows that $\cup_{u\in  \mathcal{P}_{1}\left( t_{0}\right)} u$ intersects with $t_{0}^*(\omega )$. 
To obtain the first half of the claim, it remains to compare the widths of $\cup_{u\in  \mathcal{P}_{1}\left( t_{0}\right)} u$ and $t_{0}^*(\omega )$ as both are rectangles. 
By Lemma \ref{lem:pairing}-2), one has 
\[
\textrm{Leb} \left( \pi _x\left( \bigcup_{u\in  \mathcal{P}_{1}\left( t_{0}\right)} u\right)\right) =\frac{n b_{n}b_{n+1}}{b_{n-1}}.
\]
Here, by \eqref{cond:decay3} with $k=n$, one has 
\[
\frac{n b_{n}b_{n+1}}{b_{n-1}}<\frac{b_{n+1}}{2} =\frac{nb_{n+1}}{2n}=\frac{1}{2n}\textrm{Leb} \left( \pi _x\left( t_{0} \right)\right) =\frac{1}{2} \textrm{Leb} \left( \pi _x\left( t_{0}^*(\omega) \right)\right) .
\]
Therefore, one has
\[
\pi_x\left( \bigcup_{u\in  \mathcal{P}_{1}\left( t_{0}\right)} u \right) \subset \pi_x\left( t_{0}^*(\omega )\right),
\]
and which implies the first claim. 

Since the truncation $t\mapsto t^*(\omega)$ preserves the height of $t$, by Lemma \ref{lem:pairing}-1), one has 
\[
\pi_y\left( t_{1}\right) \supset \pi_y\left( t_{0}\right)=\pi_y\left( t_{0}^*(\omega )\right).
\]
Lemma is obtained. 
\end{proof}

By the same argument as in Lemma \ref{lem:translate0}, for each $t_{1}\in \mathcal{P}_{1}(t_{0})$, one can take $\omega_{t_{1}}\in [0,(n-1)b_{n+2})$ so that every $t_{2}\in \mathcal{P}_{2}\left( t_{1}\right)$ intersects with $t_{1}^*(\omega_{t_{1}})$ in such a way that 
\begin{equation}
   \begin{cases}
   \pi_x\left( t_{2}\right) \subset \pi_x\left( t_{1}^*(\omega_{t_{1}}  )\right), \\
   \pi_y\left( t_{2}\right) \supset \pi_y\left( t_{1}^*(\omega_{t_{1}}  )\right).
   \end{cases}
\end{equation}
Note here that by Lemma \ref{lem:translate0} one has 
\[
\pi_x\left( t_{1}^*(\omega_{t_{1}})\right) \subset \pi_x\left( t_{0}^*(\omega_{t_{0}})\right)
\]
as $t_{1}^*(\omega_{t_{1}}) \subset t_{1}$, and 
\[
\pi_y\left( t_{1}^*(\omega_{t_{1}}) \right) \supset \pi_y\left( t_{0}^*(\omega_{t_{0}}) \right)
\]
as $\pi_y\left( t_{1}^*(\omega_{t_{1}}) \right) =\pi_y\left( t_{1}\right) $. 
By a recursive use of the argument above, for every chain of tiles $t_{k}\in \mathcal{P}_{k}\left(t_{k-1}\right)$, one can take $\omega_{t_{k}}\in [0,(n-1)b_{n+1+k})$ such that  
\begin{equation} \label{cond:translate0n-1}
    \begin{cases}
    \pi _x\left( t_{n-1}^*(\omega_{t_{n-1}} )\right) \subset \pi _x\left( t_{n-2}^*(\omega_{t_{n-2}} )\right) \subset \dots \subset \pi _x\left( t_{0}^*(\omega_{t_{0}} )\right) , \\
    \pi _y\left( t_{n-1}^*(\omega_{t_{n-1}} )\right) \supset \pi _y\left( t_{n-2}^*(\omega_{t_{n-2}} )\right) \supset \dots \supset \pi _y\left( t_{0}^*(\omega_{t_{0}} )\right).
\end{cases}
\end{equation}
In particular \eqref{cond:translate0n-1} yields
\begin{equation} \label{def:core}
t_{n-1}^*(\omega_{t_{n-1}} )\cap t_{n-2}^*(\omega_{t_{n-2}} )\cap \dots \cap t_{0}^*(\omega_{t_{0}} )=t_{n-1}^*(\omega_{t_{n-1}})\cap t_{0}^*(\omega_{t_{0}} ),
\end{equation}
which is a rectangle of height $b_{n}$ and width $b_{2n}$. 
In Section \ref{subsec:function}, we will define a positive function such that its support is contained in these ``core'' rectangles.

Once the translation parameters $\omega_{t_{k}}$ are chosen so that \eqref{cond:translate0n-1} is fulfilled, they will remain unchanged in the sequel.  
Henceforth, the parameter $\omega_{t_{k}}$ will be omitted from $t_{k}^*(\omega_{t_{k}} )$ and it is abbreviated as $\tau _{k}^{*}$ for simplicity of notation. 
Now, we define subsets of $[0,1)^2$ as follows. 
Given a tile $t_{0}\in \mathcal{T}_{0}$, define 
\[
B_0\left( t_{0}\right) =\tau _{0}^{*},
\]
and for each $k\in \{ 1,\dots,n-1\} $, define $B_{k}\left( t_{0}\right) \subset [0,1)^2$ by 
\begin{align*}
    B_{1}(t_{0})&=\bigcup_{t_{1}\in \mathcal{P}_{1}(t_{0})} \tau_{1}^*, \\ B_{2}(t_{0})&=\bigcup_{t_{1}\in \mathcal{P}_{1}(t_{0})} \quad \bigcup_{t_{2}\in \mathcal{P}_{2}\left(t_{1}\right)} \tau_{2}^* \\
    B_{3}(t_{0})&=\bigcup_{t_{1}\in \mathcal{P}_{1}(t_{0})} \quad \bigcup_{t_{2}\in \mathcal{P}_{2}\left(t_{1}\right)} \quad \bigcup_{t_{3}\in \mathcal{P}_{3}\left(t_{2}\right)} \tau_{3}^* \\
    &\vdots \\
    B_{n-1}(t_{0})&=\bigcup_{t_{1}\in \mathcal{P}_{1}(t_{0})} \quad \bigcup_{t_{2}\in \mathcal{P}_{2}\left(t_{1}\right)} \dots \quad \bigcup_{t_{n-2}\in \mathcal{P}_{n-2}\left(t_{n-3}\right)} \quad \bigcup_{t_{n-1}\in \mathcal{P}_{n-1}\left(t_{n-2}\right)} \tau_{n-1}^*.
\end{align*}
Set 
\[
F\left( t_{0}\right) =B_{0}\left( t_{0}\right) \cup B_{1}\left( t_{0}\right) \cup B_{2}\left( t_{0}\right) \cup \dots \cup B_{n-1}\left( t_{0}\right) ,
\]
and define $\mathcal{A}=\left\{ F\left( t_{0}\right)  \subset [0,1)^2\colon t_{0}\in \mathcal{T}_{0}\right\} $. 
For simplicity of notation, we may write $B_k$ and $F$ instead of $B_k(t_{0})$ and $F\left( t_{0}\right) $, respectively. 
\begin{figure}
\begin{center}
\includegraphics[pagebox=artbox,width=9cm]{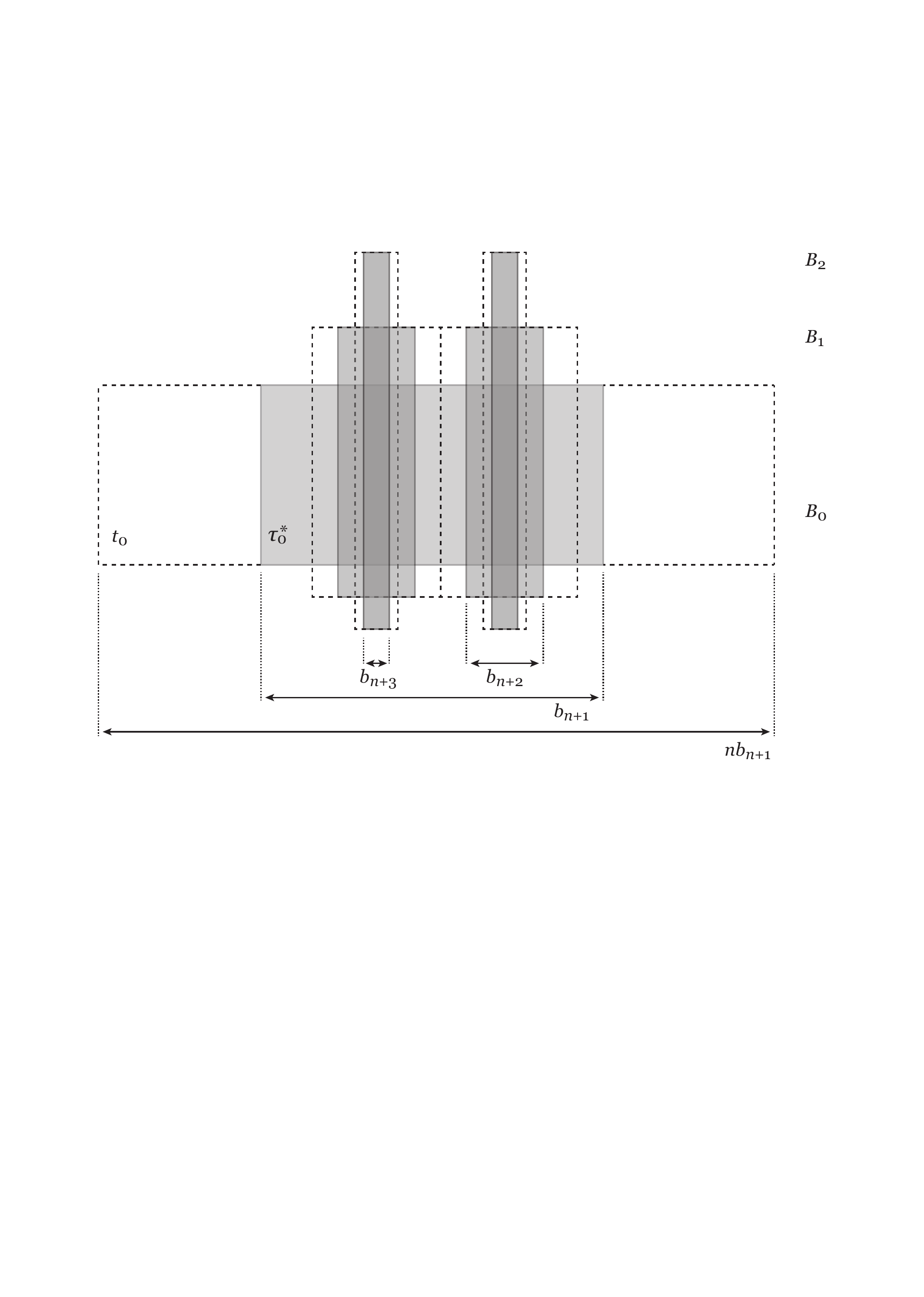}
\end{center}
\caption{Examples of $B_k$ for $k=0,1,2$.  The gray rectangles are $\tau_0^*$, $\tau_1^*$, and $\tau_2^*$. They are $1/n$'th part of the tiles (the rectangles with broken lines) $t_0$, $t_1$, and $t_2$, respectively.}
\label{fig:Bk}
\end{figure}

Let 
\[
q_k=\frac{|\tau_{k-1}^*|}{|\tau_{k}^*|}\left( =\frac{|t_{k-1}|}{|t_{k}|}=\frac{b_{n-k+1}b_{n+k}}{b_{n-k}b_{n+k+1}}\right)
\]
for $k\in \{ 1,2,\dots ,n-1\} $. 
Note that $q_k>1$ as $|\tau_{0}^*|>|\tau_{1}^*|>\dots >|\tau_{n-1}^*|$ by \eqref{cond:decay1}. 
Note also that 
\begin{equation} \label{cond:q}
    q_1\cdots q_{k-1}q_k|\tau_k^*|=q_1\cdots q_{k-1}|\tau_{k-1}^*|=\cdots =|\tau_0^*|. 
\end{equation}
By Lemma \ref{lem:pairing}-2), one sees that 
\begin{equation} \label{card_P}
    \# \mathcal{P}_{k}(t_{k-1})=\# \mathcal{P}_{k}\left(u_{k-1}\right)=q_k
\end{equation}
for every $t_{k-1},u_{k-1}\in \mathcal{P}_{k-1}\left(t_{k-2}\right)$, $k=2,\dots ,n-1$. 

In the following lemma we show that due to the big difference between their sides, the collections $B_k(t_0)$, $k=1, \dots ,n$, have very small intersection with each other. 
Hence, the total area of the figure $F$ is almost the same as the sum of individual sets $B_k(t_0)$.

\begin{lemma} \label{lem:area}
Let $F\in \mathcal{A}$ such that $F=F(t_{0})$ for some $t_{0}\in \mathcal{T}_{0}$.
\begin{enumerate}
    \item[1)] For every $k\in \{ 0,1,\dots ,n-1\}$, one has $|B_k(t_0)|=b_n b_{n+1}$. 
    \item[2)] One has 
    \[
    1-\lambda <\frac{|F|}{\sum_{k=0}^{n-1}\left| B_{k}\left( t_{0}\right) \right| }\leq 1,
    \]
    that is, the area of each $F\in \mathcal{A}$ is close to the sum of its components.
    \item[3)]
    $$
    \# \mathcal{A}\leq \frac{1}{|F|}.
    $$
\end{enumerate}
\end{lemma}
\begin{proof}
By definition and Lemma \ref{lem:pairing}-2), with the aid of \eqref{cond:q} and \eqref{card_P}, it follows that 
\[
|B_k(t_0)|=q_{1} \times q_{2} \times \dots \times q_{k}\times |\tau_k^*|=|\tau_0^*|=b_n b_{n+1}.
\]

Next, we show (2). 
Note that $\sum_{k=0}^{n-1}\left| B_{k}\left( t_{0}\right) \right| =nb_nb_{n+1}=n|\tau_0^*|$ by (1). 
By construction of $F$, one has 
\begin{align*}
    |F|
    &=|B_0|+|B_1\setminus B_0|+|B_2\setminus B_1|+\cdots +|B_{n-1}\setminus B_{n-2}| \\
    &=|\tau_0^*|+(b_{n-1}-b_n)q_1b_{n+2}+(b_{n-2}-b_{n-1})q_1q_2b_{n+3}+\cdots \\
    &\quad \cdots +(b_{1}-b_{2})q_1q_2\dots q_{n-1}b_{2n},
\end{align*}
and thus 
\begin{align*}
    \frac{|F|}{\sum_{k=0}^{n-1}\left| B_{k}\left( t_{0}\right) \right| }
    &=\frac{|F|}{n|\tau_0^*|} \\
    &=\frac{1}{n} \Big\{ 1+\frac{(b_{n-1}-b_n)q_1b_{n+2}}{|\tau_0^*|}+\frac{(b_{n-2}-b_{n-1})q_1q_2b_{n+3}}{|\tau_0^*|} +\cdots \\
    &\quad \cdots +\frac{(b_{1}-b_{2})q_1q_2\cdots q_{n-1}b_{2n}}{|\tau_0^*|} \Big\}.
\end{align*}
Here, for every $k\in \{ 1,\dots ,n-1\}$, it follows from \eqref{cond:q} and \eqref{cond:decay3} that 
\begin{align*}
    \frac{(b_{n-k}-b_{n-k+1})q_1\cdots q_kb_{n+k+1}}{|\tau_0^*|} &=\frac{(b_{n-k}-b_{n-k+1})q_1\cdots q_kb_{n+k+1}}{q_1\cdots q_k|\tau_k^*|} \\ &=\frac{b_{n-k}b_{n+k+1}}{|\tau_k^*|} -\frac{b_{n-k+1}b_{n+k+1}}{|\tau_k^*|} \\
    &=1-\frac{b_{n-k+1}}{b_{n-k}} >1-\lambda .
\end{align*}
Hence 
\[
\frac{|F|}{\sum_{k=0}^{n-1}\left| B_{k}\left( t_{0}\right) \right| }=\frac{|F|}{n|\tau_0^*|} >\frac{n-(n-1)\lambda}{n} >1-\lambda .
\]
The opposite estimate $|F|\leq \sum_{k=0}^{n-1}\left| B_{k}\left( t_{0}\right) \right|$ is clear. 

By the definition of the set $\mathcal{T}_{0}$ and (2) it follows that
\[
\# \mathcal{A}=\#\{F(t_0)\colon t_0 \in \mathcal{T}_{0}\}=\#\mathcal{T}_{0}=\frac{1}{nb_n b_{n+1}}\leq \frac{1}{|F|}.
\]
Lemma is obtained. 
\end{proof}

\begin{lemma} \label{lem:figure_enlarged}
%Given $\varepsilon \in (0,1)$, one can choose $\{ b_k \}_{k=1}^{2n}$ such that 
One has 
\[
\left| \bigcup_{F\in \mathcal{A}}F\right| \geq \frac{1}{2}.
\]
\end{lemma}
\begin{proof}
Since every tile $t$ participates in the pairing procedure described above and each tile is paired with exactly one set $F(t_{0})$, then
$$
\bigcup_{F\in \mathcal{A}}F =\bigcup_{k=0}^{n-1} \left\{ \tau^*(=t^*(\omega_{t}))\colon t \in \mathcal{T}_{k}\right\}.
$$
But by Lemma \ref{Lmm-cov} we have that the last set has measure larger than $1/2$.
\end{proof}

\section{Main Proposition} \label{sec:mainprop}

We are now ready to prove a key Proposition:

\begin{proposition} \label{prop:main1}
For every $\varepsilon >0$ and every $n>0$ which is a power of $2$, there exists a function $f\in L^\infty([0,1)^2)$, with $f\geq 0$ and $\|f\|_{L^1}\leq 1+o_{\varepsilon}(1)$, for which there exists two disjoint subsets $E,N\subset [0,1)^2$, with $|E|\geq 1/3$ and $|N|\leq \varepsilon$, such that for every $z\in E$ there exists a dyadic rectangle $R$ from $\mathcal{F}_D$, with $z \in R$ such that
\begin{equation} \label{cond:div0}
    \frac{1}{|R|}\int_{R} f\, dxdy \geq \frac{n}{2},
\end{equation}
and for every $z\in [0,1)^2 \setminus N$ and any dyadic rectangle $R \in \mathcal{F}_C$, with $z\in R$, we have that 
\begin{equation} \label{cond:conv0}
    \frac{1}{|R|}\int_{R} f\, dxdy \leq 3.
\end{equation}
\end{proposition}

%Here we assume that decay in \eqref{cond:decay1} is fast enough so that Lemma \ref{lem:difference} is satisfied. 
%Recall that each $F=F\left( t_{0}\right) \in \mathcal{A}$, $t_{0}\in \mathcal{T}_0$, is determined by a chain of tiles from $\mathcal{P}_{1}\left(t_{0}\right),\mathcal{P}_{2}\left(t_{1}\right),\dots ,\mathcal{P}_{n-1}\left(t_{n-2}\right)$, where each $t_{k}$ is taken from $\mathcal{P}_{k}\left(t_{k-1}\right) \subset \mathcal{T}_{k}$ for $k\in \{ 1,\dots ,n-2\}$. 
%One can find certain geometric structure on $F$ due to \eqref{cond:decay3}. 
In Section \ref{subsec:function}, we will define a positive function associated with $F\in \mathcal{A}$ using the geometric structure \eqref{cond:translate0n-1}. 
Then we will prove some preliminary estimates in Section \ref{subsec:preliminary}, and prove Proposition \ref{prop:main1} in Section \ref{subsec:proofmain_prop}. 

\subsection{Geometric construction of a function} \label{subsec:function}
\begin{figure}
\begin{center}
\begin{subfigure}[t]{0.45\textwidth}
\includegraphics[pagebox=artbox,width=\textwidth]{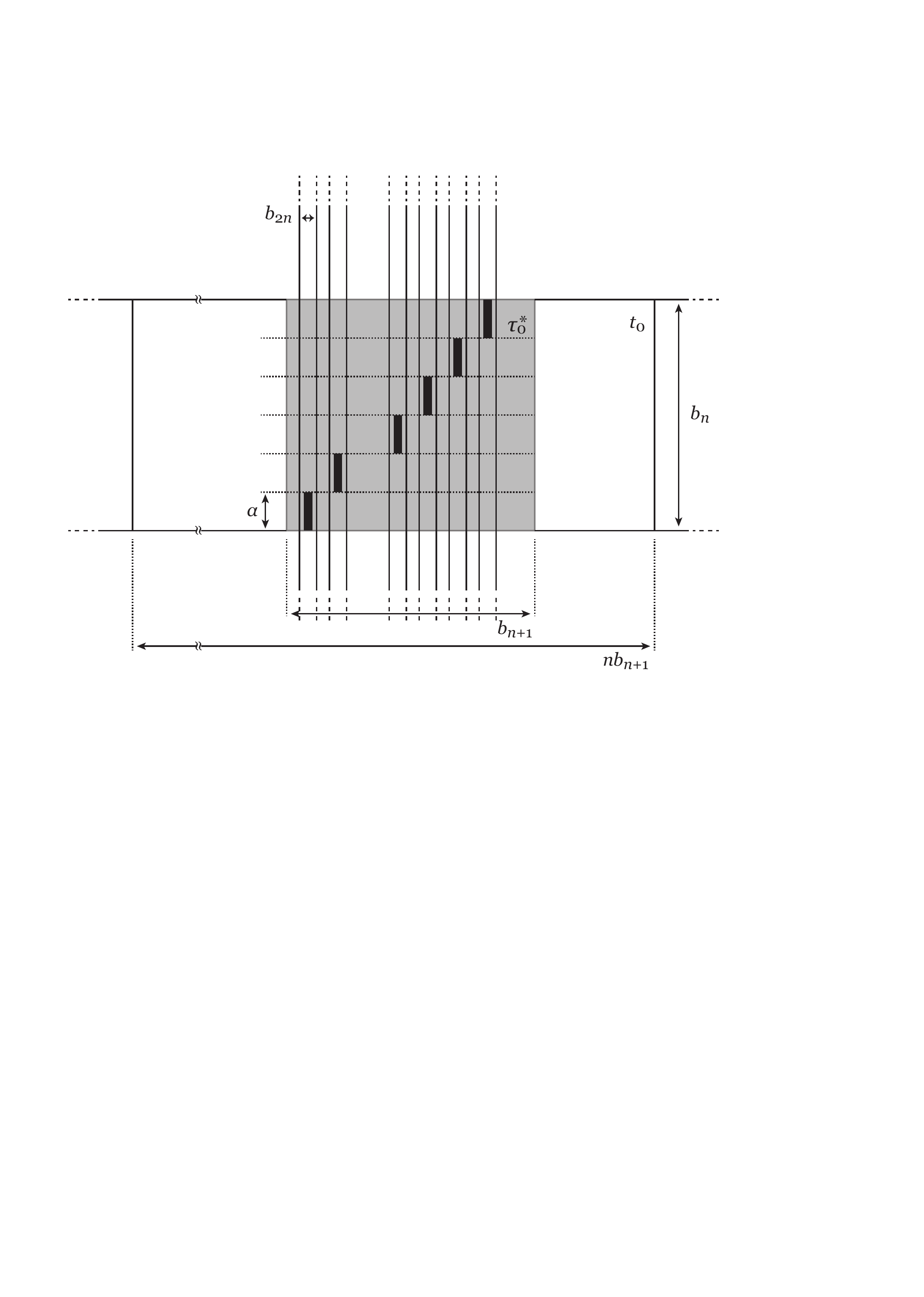}
\caption{The gray rectangle is $\tau _0^*$, and the vertical, long rectangles are $\tau_{n-1}^*$. Each vertical one has sides $b_1$ (height) and $b_{2n}$ (width). $\alpha =b_1b_{2n}/b_{n+1}$. }
\label{fig:support_local}
\end{subfigure}
\hfill
\begin{subfigure}[t]{0.45\textwidth}
\includegraphics[pagebox=artbox,width=\textwidth]{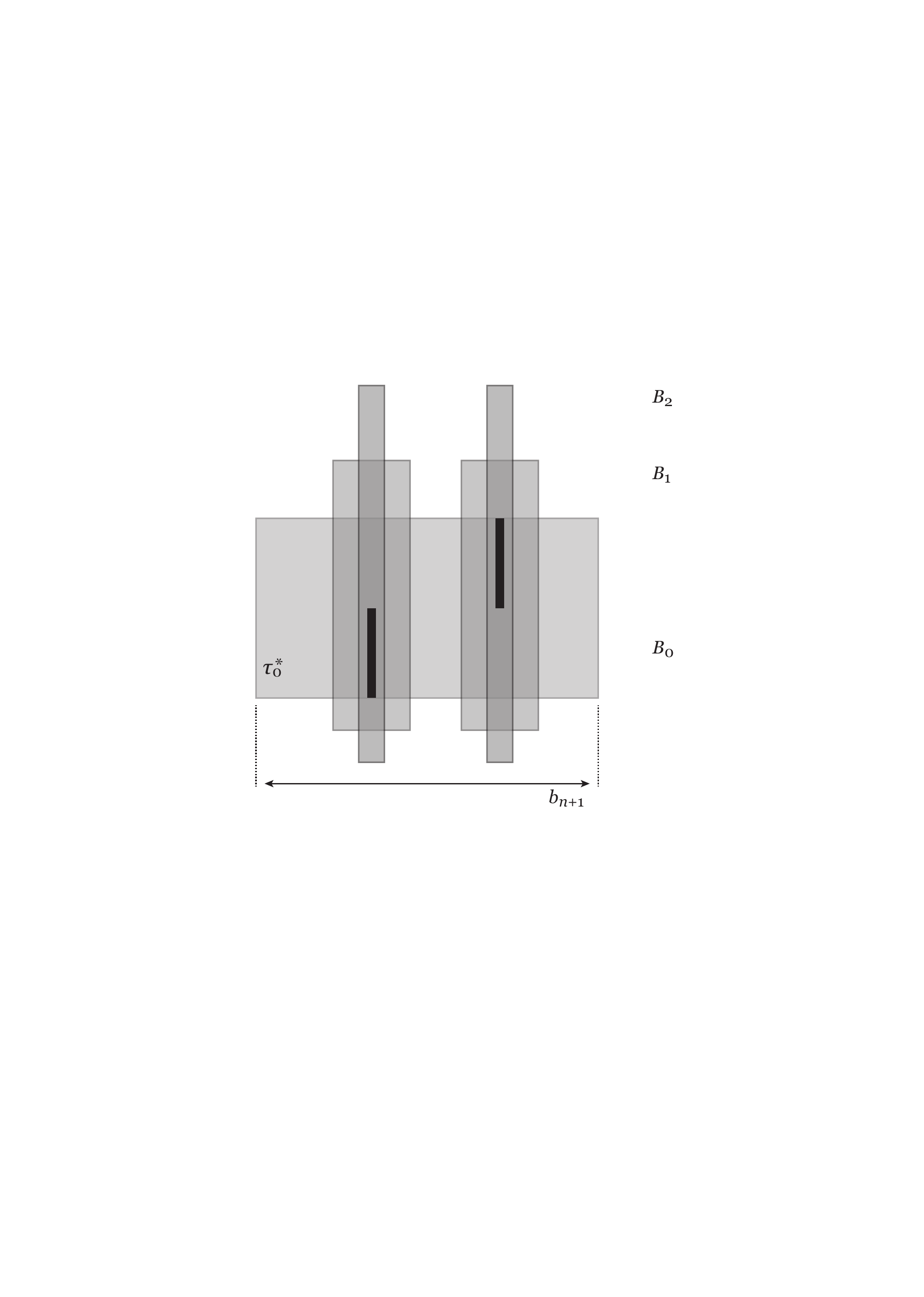}
\caption{Example with $3$ generations}
\label{fig:support_global}
\end{subfigure}
\caption{Placement of $s_\rho$ (black rectangles). }
\label{fig:support}
\end{center}
\end{figure}

Let $\{ b_k\} _{k=1}^{2n}$ be a decreasing sequence satisfying \eqref{cond:decay1}, \eqref{cond:decay2} with $\lambda <\frac{1}{2n(n-1)}$, and \eqref{cond:decay3}. 
Recall that each $F=F\left( t_{0}\right) \in \mathcal{A}$, $t_{0}\in \mathcal{T}_0$, is determined by a chain of tiles from $\mathcal{P}_{1}\left(t_{0}\right),\mathcal{P}_{2}\left(t_{1}\right),\dots ,\mathcal{P}_{n-1}\left(t_{n-2}\right)$, where each $t_{k}$ is taken from $\mathcal{P}_{k}\left(t_{k-1}\right) \subset \mathcal{T}_{k}$ for $k\in \{ 1,\dots ,n-2\}$. 
For each $t_{0}\in \mathcal{T}_0=\mathcal{T}_{b_n,n b_{n+1}}$, recall that there are $q_1\cdots q_{n-1}$ many ``core" rectangles $\tau_{n-1}^*\cap \tau_{0}^*$.
See \eqref{def:core} and \eqref{card_P}. 
To obtain an aimed function associated with $F=F(t_0)\in \mathcal{A}$, for each $t_0 \in \mathcal{T}_0$, we will place a mass at each core rectangle in such a way that they are distributed uniformly along the vertical direction. 
More specifically, we do as follows. 
Let $t_0 \in \mathcal{T}_0$, and let $c_n=q_1\cdots q_{n-1}$. 
Here recall that 
\[
q_k=\frac{|\tau_{k-1}^*|}{|\tau_{k}^*|}\left( =\frac{|t_{k-1}|}{|t_{k}|}=\frac{b_{n-k+1}b_{n+k}}{b_{n-k}b_{n+k+1}}\right)
\]
for $k\in \{ 1,2,\dots ,n-1\} $, and hence
\[
c_n=q_1\cdots q_{n-1}=\frac{|\tau_{0}^*|}{|\tau_{n-1}^*|} =\frac{b_{n}b_{n+1}}{b_{1}b_{2n}}
\]
as we have seen in \eqref{cond:q}. Note that $c_n$ is the total number of rectangles $\tau_{n-1}^*$ that intersect with $\tau_0^*$. 
%We have $q_k>1$ and \eqref{card_P}. 
We enumerate core rectangles as 
\[
\{\tau_{n-1}^*\cap \tau_{0}^*\colon t_{n-1}\in \mathcal{P}_{n-1}\left(t_{n-2}\right)\} =\{ \gamma_j^*\colon j=0,1,\dots ,c_n-1\}. 
\]
For $\rho \in (0,b_{2n})$, define the set $s_\rho $ as the rectangle with side length $\rho $ (width) and $\alpha =b_1b_{2n}/b_{n+1}$ (height). 
For each $\gamma_j^*$, we place a single $s_\rho$ inside it in such way that
\begin{equation} \label{cond:support_v}
    \begin{cases}
    \pi _x(s_\rho)\subset \pi _x(\gamma_j^*), \\[2mm]
    \pi _y(s_\rho)=\pi _y\left(\left[ j\dfrac{b_n}{c_n},(j+1)\dfrac{b_n}{c_n}\right)\right) . 
   \end{cases}
\end{equation}
The exact $x$ coordinate of $s_\rho$ inside $\gamma_j^*$ is not important. 
Since 
\[
\frac{b_n}{c_n} =\frac{b_{1}b_{2n}}{b_{n+1}} =\alpha ,
\]
one can achieve \eqref{cond:support_v}, and hence the rectangles are distributed uniformly along vertical direction such that 
\[
\textrm{Leb} \left(\pi_y\left( \bigcup _{j=1}^{c_n-1}s_{\rho}\right) \right) =b_n. 
\]
See Figure \ref{fig:support}. 
%\begin{equation} \label{cond:support}
%    \begin{cases}
%    \pi _x(s_\rho)\subset \pi _x(\tau_{n-1}^*), \\
%    \pi _y(s_\rho)\subset \pi_y(\tau_0^*)(=\pi_y(t_0)). 
%   \end{cases}
%\end{equation}
%Hence, each $\tau_{n-1}^*$ contains exactly one $s_\rho$, and thus by \eqref{card_P} there will be $q_1\cdots q_{n-1}$ many squares $s_\rho$ among each $\tau _0^*\subset F(t_0)$. 
Now, for each $t_0 \in \mathcal{T}_0$ and the associated $F=F(t_0)\in \mathcal{A}$, we define a positive function $h=h_{t_0}$ by
\begin{equation} \label{def:function_h}
h_{t_0}(z)=\begin{cases} \dfrac{n b_{n} b_{n+1}}{S q_1\cdots q_{n-1}}, &\text{if } z\in s_\rho \subset \gamma _j^*\text{ for some } \gamma _j^*, \\
0, &\text{otherwise}, \end{cases}
\end{equation}
where $S=|s_\rho|=\rho \alpha$, the area of $s_\rho$. The purpose of distributing the support of $h_{t_0}$ like this will be clear in Lemma \ref{lem:concentration_support}, (2) below.

Hence $\text{supp}h\subset \tau _0^*\subset t_0$. 
Recall that by Lemma \ref{lem:area} we have for every $F\in \mathcal{A}$ that 
\begin{equation} \label{cond:areaF}
    n b_{n} b_{n+1}\geq |F|\geq (1-\lambda)n b_{n} b_{n+1}.
\end{equation}
Thus $\int_F h_{t_0}\,dxdy$ is close to $|F|$.
Since each $F=F(t_0)\in \mathcal{A}$ is indexed by  $t_0\in \mathcal{T}_0$, the definition of function $h_{t_0}$ in \eqref{def:function_h} can be extended to a positive function $f$ on the unit square as follows
\begin{equation} \label{def:function_f}
f(z)=\sum_{t_0 \in \mathcal{T}_0}h_{t_0}(z).
\end{equation}

\begin{lemma} \label{lem:concentration_support}
Let $f$ be defined as in \eqref{def:function_f}. 
\begin{enumerate}
    \item[1)] For every $t_0 \in \mathcal{T}_0$, $k\in \{ 0,1,\dots ,n-1\} $ and $t\in \mathcal{P}_{k}\left(t_{k-1}\right)$, with a convention that $t=t_{0}$ for $k=0$, we have 
    \[
    \frac{1}{|\tau^*|}\int_{\tau^*} f\, dxdy =n,
    \]
    and 
    \[
    \frac{1}{|t|}\int_{t} f\, dxdy =1.
    \]
    \item[2)] For every $t_0 \in \mathcal{T}_0$ and a dyadic rectangle $A\subset t_0$ with sides $x,y$, so that $x\leq b_n$, $y = n b_{n+1}$ we have that
    \[
    \frac{1}{|A|}\int_A f\, dx dy=1.
    \]
\end{enumerate}
\end{lemma}
\begin{proof}
By construction, one has 
\begin{align*}
    \int_{\tau_{k}^*}f\, dxdy=\int_{\tau_{k}^*}h\, dxdy 
    &=\frac{n b_{n} b_{n+1}}{S q_1\cdots q_{n-1}}\times \left| \text{supp}h\cap \tau_{k}^*\right| \\
    &=\frac{n b_{n} b_{n+1}}{S q_1\cdots q_{n-1}}\times S \times q_{k+1}\times \dots \times q_{n-1} \\
    &=\frac{n b_{n} b_{n+1}}{q_1\cdots q_{k}}.
\end{align*}
Thus we have
\[
\frac{1}{|\tau_{k}^*|} \int_{\tau_{k}^*}f\, dxdy= \frac{1}{|\tau_{k}^*|} \frac{n b_{n} b_{n+1}}{q_1\cdots q_{k}} = \frac{n b_{n} b_{n+1}}{b_{n} b_{n+1}} =n 
\]
as $q_1\cdots q_{k-1}q_{k}|\tau_{k}^*|=|\tau_{0}^*|=b_{n} b_{n+1}$ from \eqref{cond:q}. 
%Similarly, it follows from \eqref{card_P} that
%\[
%\frac{1}{|t_{(k)}^*|} \int_{t_{(k)}^*}h\, d m \geq \frac{1}{|t_{(k)}^*|} \frac{n b_{n} b_{n+1}}{q_1\cdots q_{k}} \left( 1-\frac{1}{q_{k+1}}\right)\cdots \left( 1-\frac{1}{q_{n-1}}\right) = n\prod _{j=k+1}^{n-1}\left( 1-\frac{1}{q_{j}}\right) .
%\]
%

For the latter half of (1), note that 
\[
\int_{t} f\, dxdy=\int_{t} h\, dxdy=\int_{\tau^*} h\, dxdy=\int_{\tau^*} f\, dxdy  
\]
since $\text{supp}h$ is contained in the core rectangles $s_\rho \subset \tau_{n-1}^*\cap \tau_{0}^*\subset \tau^*$ by \eqref{cond:support_v}. 
It then follows from the first assertion that
\[
\int_{t} f\, dxdy=\int_{\tau^*} f\, dxdy=n|\tau ^*|=|t|. 
\]
%and 
%\[
%\int_{t} h\, d m=\int_{t^*} h\, d m\geq n|t^*|\prod _{j=k+1}^{n-1}\left( 1-\frac{1}{q_{j}}\right) =|t| \prod _{j=k+1}^{n-1}\left( 1-\frac{1}{q_{j}}\right).
%\]

Next, we show (2). 
Since the support of $f$ is (or the rectangles $s_\rho$) distributed uniformly along the vertical direction, it follows that 
\[
\int_A f\, dx dy=\frac{x}{b_n}\int_{t_0} f\, dx dy.
\]
Thus, in view of (1) of this lemma, one has
\[
\frac{1}{|A|}\int_A f\, dx dy=\frac{1}{|A|} \frac{x}{b_n}\int_{t_0} f\, dx dy=\frac{1}{|A|} \frac{x}{b_n} |t_0|=\frac{x n b_n b_{n+1}}{xyb_n} =\frac{n b_{n+1}}{y} =1.
\]
Lemma is obtained. 
\end{proof}

\begin{comment}
   We now determine the positions of $\text{supp}h$, namely those of the squares $s_\rho \subset \tau_{n-1}^*\cap \tau_{0}^*$. 
Observe that each $F=F\left( t_{0}\right) \in \mathcal{A}$, $t_{0}\in \mathcal{T}_0$, one sees that 
\[
\# \left\{ t_{n-1}\subset F \colon t_{n-1}\in \mathcal{T}_{n-1} \right\} =\# \mathcal{P}_{1}\left(t_{0}\right)\times \# \mathcal{P}_{2}\left(t_{1}\right)\times\dots \times \# \mathcal{P}_{n-1}\left(t_{n-2}\right),
\]
where note that $\# \mathcal{P}_{k}\left(t_{k-1}\right)=\# \mathcal{P}_{k}\left(u_{k-1}\right)$ for every $t_{k-1},u_{k-1}\in \mathcal{P}_{k-1}\left(t_{k-2}\right)$, $k=2,\dots ,n-1$.  
We set $c_n=c_n(t_{0})=\# \mathcal{P}_{1}\left(t_{0}\right)\times \# \mathcal{P}_{2}\left(t_{1}\right)\times\dots \times \# \mathcal{P}_{n-1}\left(t_{n-2}\right)$, and denote by
\[
\left\{ t_{n-1}\subset F \colon t_{n-1}\in \mathcal{T}_{n-1} \right\} =\left\{ t_{n-1,j}\colon j=0,1,\dots ,c_n-1\right\} .
\]
%from left to right. 
One can assume $\rho \in (0,\min \{ b_{2n},b_{n}/c_{n}\})$. 
Then for each $\tau_{n-1,j}^*$, $j\in \{0,1,\dots ,c_n-1\}$, the square $s_\rho \subset \tau_{n-1,j}^*\cap \tau_{0}^*$ is placed in such a way that 
\begin{equation} \label{cond:support}
    \begin{cases}
    \pi _x(s_\rho)\subset \pi _x(\tau_{n-1,j}^*) \text{ and } \\
    \pi _y(s_\rho)\subset \left[ j\dfrac{b_n}{c_n},(j+1)\dfrac{b_n}{c_n}\right). 
   \end{cases}
\end{equation}
Hence the squares $s_\rho$ are distributed uniformly along vertical direction. 
%Let \[s_n=\frac{b_n}{\underline{a} _{n+1}} .\]
\end{comment}

\subsection{Preliminary estimates} \label{subsec:preliminary}

Recall the definitions \eqref{def:C} and \eqref{def:D} and define the collections
\begin{equation}\label{c1}
\mathcal{C}=\left\{\frac{1}{2^k}\in [0,1]\colon k \in C\right\}
\end{equation}
and
\begin{equation}\label{d1}
\mathcal{D}=\left\{\frac{1}{2^s}\in [0,1]\colon s \in D\right\}.
\end{equation}

Let $\mathcal{D}=\{ b_n\}_{n\in \mathbb{N}}$. For each $b_n \in \mathcal{D}$, we define
\begin{align*}
\overline{a} _n&=\sup \left\{ a\in \mathcal{C}\colon a<b_n\right\} , \\
\underline{a} _n&=\inf \left\{ a\in \mathcal{C}\colon a>b_n\right\}.
\end{align*}

Note that condition \eqref{cond:dst} is equivalent to the following condition 
\begin{equation} \label{regularity}
\liminf_{n \rightarrow \infty}\left( \max \left\{ \frac{\bar a_n}{b_n},\frac{b_n}{\underline{a}_n}\right\} \right) =0.
\end{equation}
We have the following lemma. 
The proof is a direct consequence of \eqref{regularity}, and is omitted.

\begin{lemma} \label{lem:distortion}
Assume the condition \eqref{regularity}. 
Given $b_{1}>\dots >b_{n}$, and $\lambda \in (0,1)$, one can choose $b_{n+1}>\dots >b_{2n}$ such that for every $k=0,1,\dots ,n-1$
\begin{enumerate}
    \item[1)] $\dfrac{\overline{a} _{n+1}}{b_{n+1}} \leq \lambda b_{n}$, \hbox{ and }$\dfrac{b_{n+k}}{\underline{a}_{n+k}}\leq \dfrac{1}{n}$
\end{enumerate}

\begin{enumerate}
    \item[2)] $\dfrac{\overline{a} _{n+k+1}}{b_{n+k+1}} \leq \lambda \dfrac{b_{n-k}}{\overline{a} _{n-k-1}}$, and $\dfrac{b_{n+k+1}}{\underline{a} _{n+k+1}} \leq \lambda \dfrac{\overline{a} _{n-k}}{b_{n-k}} $,
\end{enumerate}
where $\overline{a} _{0}=1$ as a convention. 
\end{lemma}

We denote the axis-parallel \emph{dyadic} rectangles with side lengths $x$ and $y$ by $A_{x, y}$, where $x$ is the length of vertical side and $y$ is that of horizontal one. 
Let 
\[
\mathcal{W}_{C}=\left\{ A_{x,y}\in \mathcal{F}_C\colon A_{x,y}\cap \text{supp} f\neq \emptyset \right\} ,
\]
where $f$ is the positive function defined in \eqref{def:function_f}. 

\begin{lemma} \label{lem:def_N} \label{lem:removalfamily}
Assume we have \eqref{regularity}. 
Let $f$ be a positive function defined as in \eqref{def:function_f}. 
Given $\varepsilon \in (0,1)$, one can choose a sufficiently small $\lambda \in (0,1)$ and a set $N\subset [0,1)^2$ in such a way that the following properties hold: 
\begin{enumerate}
    \item[1)] For every $z\in [0,1)^2\setminus N$ and $A\in \mathcal{W}_{C}$ with $A\ni z$ one has 
    \[
    \frac{1}{|A|} \int _{A} f\, dxdy\leq 3.
    \]
    \item[2)] $|N|<\varepsilon$. 
\end{enumerate}
\end{lemma}
\begin{proof}
We divide our argument into the following two cases with respect to the range of $x$: 1) $x\in (b_n,1]$ and 2) $x\leq b_{n}$. 
Depending on the case we will define sets denoted by $N_1,N_2\subset [0,1)^2$.
For these sets the integral averages are expected to be large so they will constitute $N$ and will need to be removed.
\begin{enumerate}
    \item[Case 1)] $x\in (b_n,1]$: Then there is $k\in \{ 1,2,\dots ,n-1,n\}$ such that $x\in (b_{n-k+1},b_{n-k}]$, with a convention $b_0=1$. 
    \begin{enumerate}
        \item[1-i)] Assume $y>n b_{n+k}$: 
        By Lemma \ref{lem:concentration_support}-1), for each rectangle $A=A_{x,y}$, one has
        \[
        \int _{A} f(z)\, d z\leq \sum_{\substack{t\in \mathcal{T}_{b_{n-k+1}, n b_{n+k}} \colon \\ t\cap A \neq \emptyset }}\int _{t} f(z)\, d z\leq \sum_{\substack{t\in \mathcal{T}_{b_{n-k+1}, n b_{n+k}} \colon \\ t\cap A \neq \emptyset }} |t|.
        \]
        Here observe that 
        \begin{equation} \label{ineq:card}
            \# \left\{ t\in \mathcal{T}_{b_{n-k+1}, n b_{n+k}} \colon t\cap A \neq \emptyset \right\} \leq 3\frac{|A|}{|t|} .
        \end{equation}
\begin{figure}[H]
\begin{center}
\includegraphics[pagebox=artbox,width=6cm]{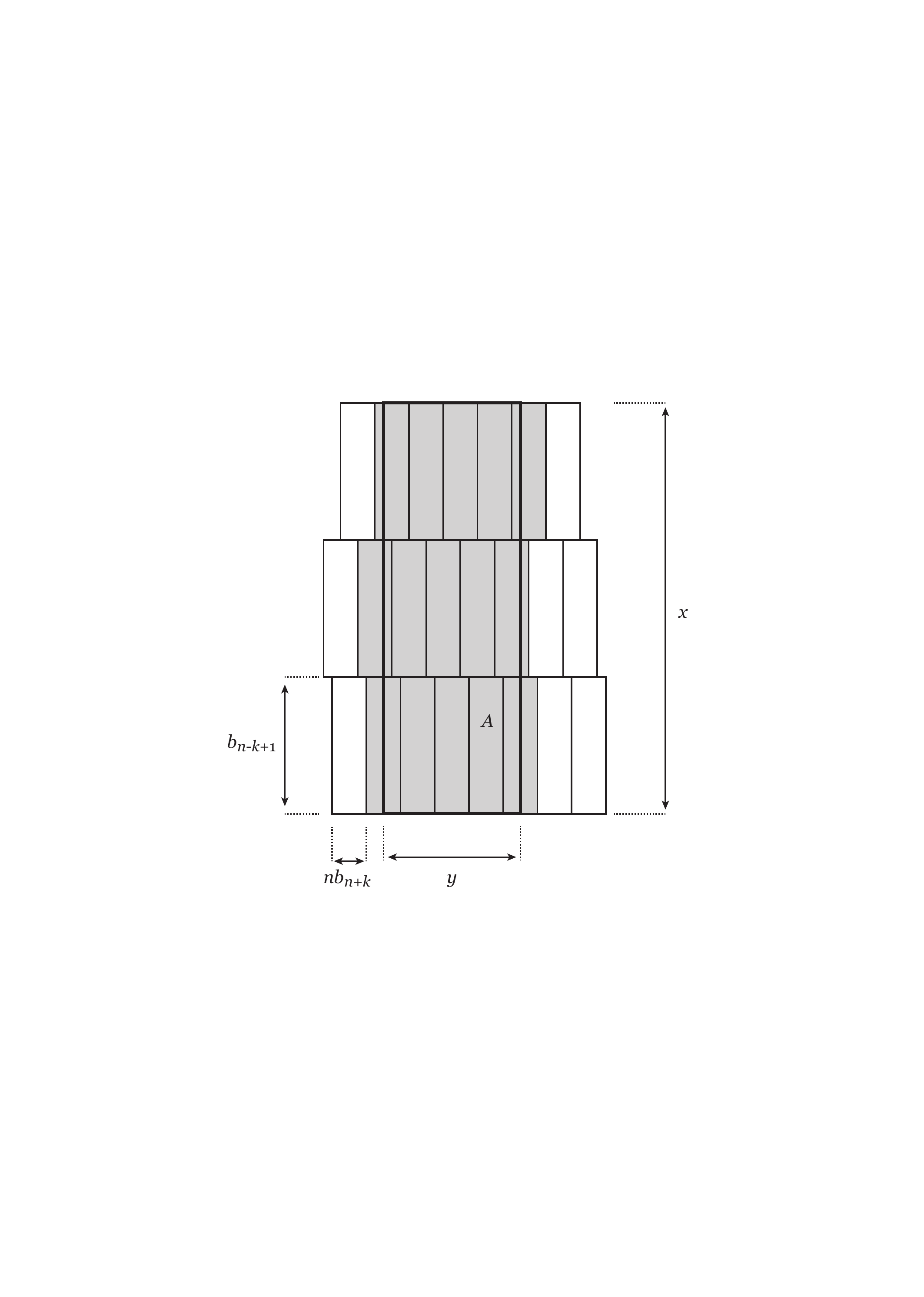}
\end{center}
\caption{Intuition behind formula \eqref{ineq:card}: one needs to estimate the number of gray rectangles that intersect $A=A_{xy}$ partially. Since $b_{n-k+1}<x$ and $n b_{n+k}<y$ the total measure of such rectangles is small.}
\label{fig:ubound}
\end{figure}
        Indeed, since $\pi_y(A)$ is dyadic, one has 
        \begin{align*}
            \# \left\{ t\in \mathcal{T}_{b_{n-k+1}, n b_{n+k}} \colon t\cap A \neq \emptyset \right\} 
            &\leq { \frac{x}{b_{n-k+1}} \left(  \frac{y}{n b_{n+k}} +2\right)} \\
            &=\frac{xy}{n b_{n-k+1} b_{n+k}} +2\frac{x}{b_{n-k+1}} \\
            &<\frac{xy}{n b_{n-k} b_{n+k+1}} +2\frac{x}{b_{n-k+1}} \frac{y}{n b_{n+k}} =3\frac{|A|}{|t|}.
        \end{align*}
        Thus, it follows that 
        \[
        \frac{1}{|A|} \int _{A} f(z)\, d z\leq \frac{1}{|A|} \times 3\frac{|A|}{|t|} \times |t| =3. 
        \]
        
        \item[1-ii)] Assume $y\leq n b_{n+k}$. Then, by Lemma \ref{lem:distortion}, (1) second inequality, we cannot have $b_{n+k}<y\leq nb_{n+k}$ (recall that $y\in \mathcal{C}$). Hence, we can assume that $y< b_{n+k}$: 
        For $k\in \{1,2,\dots ,n\}$, define 
        \[
        N_{1,k}=\bigcup_{\substack{A_{x,y}\in \mathcal{W}_{C}\colon \\ x\in (b_{n-k+1},b_{n-k}],\hbox{ } y\in (0,b_{n+k}]}} A_{x,y},
        \]
        and 
        \[
        N_1=\bigcup _{k=1}^{n} N_{1,k}.
        \]
        (Recall that we let $b_0=1$.)
    \end{enumerate}
    \item[Case 2)] $x\leq b_{n}$: 
    \begin{enumerate}
        \item[2-i)] $y>n b_{n+1}$: 
        For each rectangle $A=A_{x,y}$ due to the dyadicity of all rectangles involved
        \[
        \frac{1}{|A|} \int_{A}f(z)\,dz =\frac{1}{|A|} \sum_{\substack{t_0 \in \mathcal{T}_0;\\ t_0\cap A\neq \emptyset}}\int_{A\cap t_0}f(z)\,dz = \frac{1}{|A|} \frac{y}{nb_{n+1}} \int_{t_0\cap A}f(z)\,dz 
        \]
        as $\# \{ t_0\in \mathcal{T}_0\colon t_0\cap A\neq \emptyset \}=y/(nb_{n+1})$. %, where $t_0\in \mathcal{T}_0$ is a tile that intersects $A$.  
        For $t_0\cap A(\neq \emptyset)$, by Lemma \ref{lem:concentration_support}-2) one has 
        \[
        \int_{t_0\cap A}f(z)\,dz=|t_0\cap A|=\frac{x}{b_n} |t_0|
        \]
        since the support of $f$ (or the rectangles $s_\rho $)  is distributed uniformly along the vertical direction.  
        It follows that 
        \[
        \frac{1}{|A|} \int_{A}f(z)\,dz = \frac{1}{|A|} \frac{y}{nb_{n+1}} \frac{x}{b_n} |t_0|=\frac{|t_0|}{n b_n b_{n+1}}=\frac{n b_n b_{n+1}}{n b_n b_{n+1}}=1.
        \]
        \item[2-ii)] $y\leq n b_{n+1}$. Similar to above  by Lemma \ref{lem:distortion}, we can assume that $y< b_{n+1}$: Define
        \[
        N_2=\bigcup_{\substack{A_{x,y}\in \mathcal{W}_{C}\colon \\ x\in (0,b_{n}],\hbox{ }y \in (0,b_{n+1}] }} A_{x,y}. 
        \]
        Note in fact that we have $N_2\subset N_1$. Indeed, let $x\in (0,b_n]$ $y\in (0,b_{n+1}]$. 
        Then there is $k\in \{1,\dots ,n-1\}$ such that $y\in (b_{n+k+1},b_{n+k}]$, with the convention that $b_{2n+1}=0$. 
        For each $y\in (b_{n+k+1}, b_{n+k}]$, one has $A_{x,y}\subset A_{\tilde{x},y}$ for any $\tilde{x}=\tilde{x}_k\in (b_{n-k+1},b_{n-k}]$ since $x\leq b_n\leq b_{n-k+1}<\tilde{x}_k$. 
        Since $A_{\tilde{x}_k,y}$ belongs to $N_{1,k}$ defined in case 1-ii), one has  $N_2\subset N_1$. 
    \end{enumerate}
\end{enumerate}
It follows that letting 
\[
N=N_1\cup N_2
\]
will imply the first assertion. 
Note that $N=N_1$ since $N_2\subset N_1$ as observed. 

Next, we will show (2) which claims that the Lebesgue measure of $N=N_1$ can be made arbitrarily small. 
To see it, given $F\in \mathcal{A}$, we set $N(F)=\cup _{\{A\in N\colon A\cap F\neq \emptyset\}}A$, and show $|N(F)|<\varepsilon |F|$. 
Once it is shown, we have (2) since $\# \mathcal{A}\leq 1/|F|$ by Lemma \ref{lem:area}-3), and hence $|N|\leq \#\mathcal{A}\times |N(F)|<\varepsilon$. 
Below, we let $N_{1,k}(F)=\cup _{\{A\in N_{1,k}\colon A\cap F\neq \emptyset\}}A$ for $k\in \{ 1,\dots ,n\}$. 
Thus $N(F)=\cup _{k=1}^{n}N_{1,k}(F)$ as $N=N_1$. 

For each $k\in\{ 1,\dots ,n\}$, the rectangle $A=A_{x,y}$, with $x\in(b_{n-k+1},b_{n-k}]$, $y\in (0, b_{n+k}]$ satisfies 
\begin{equation} \label{cond:areaA}
    |A_{x,y}|\leq xy\leq \overline{a} _{n-k}\overline{a} _{n+k} \leq \lambda b_{n-k+1} b_{n+k}
\end{equation}
by Lemma \ref{lem:distortion}-2). 
Note here that $b_{n-k+1} b_{n+k}$ is the area of $\tau_{k-1}^*$. 
Hence, repeating the same argument as in the proof of Lemma \ref{lem:area}, by replacing $B_{k-1}$ with $N_{1,k}(F)$, we will get that 
\[
|N(F)|=\left| \bigcup _{k=1}^{n}N_{1,k}(F)\right| \leq \sum _{k=1}^{n}|N_{1,k}(F)|.  
\]
Here for each $k\in \{1,\dots,n\}$ one has
\[
|N_{1,k}(F)|\leq q_1\cdots q_{k-1}\overline{a} _{n-k}\overline{a} _{n+k} \leq \lambda q_1\cdots q_{k-1} b_{n-k+1} b_{n+k}=\lambda b_{n} b_{n+1}
\]
by \eqref{cond:areaA} and \eqref{cond:q}. 
It follows that 
\[
|N(F)| \leq \sum _{k=1}^{n}|N_{1,k}(F)|\leq \lambda \sum _{k=1}^n b_{n} b_{n+1} =\lambda b_{n} b_{n+1}. 
\]
Thus by \eqref{cond:areaF}, or Lemma \ref{lem:area}, one has
\[
|N(F)|\leq \lambda b_{n} b_{n+1} \leq \frac{\lambda}{1-\lambda} |F|.
\]
Lemma is obtained. 
\end{proof}

\subsection{Proof of Proposition \ref{prop:main1}} \label{subsec:proofmain_prop} 

Let $f$ be the positive function defined in \eqref{def:function_f}. 
By construction and Lemma 6-1) it follows that for every $\tau^*$ we have that
\[
\frac{1}{|\tau^*|}\int_{\tau^*}f\, dxdy =n.
\]
However, the rectangle $\tau^*$ is not dyadic since neither is $\pi_x(\tau^*)$. 
Note, however, that $\pi_y(\tau^*)$ is dyadic. 
In order to fix this issue we now consider two dyadic, adjacent rectangles that are horizontal translations of $\tau^*$ and that cover $\tau^*$. 
See Figure \ref{fig:dyadicLR} below.  
\begin{figure}[H]
\begin{center}
\includegraphics[pagebox=artbox,width=9cm]{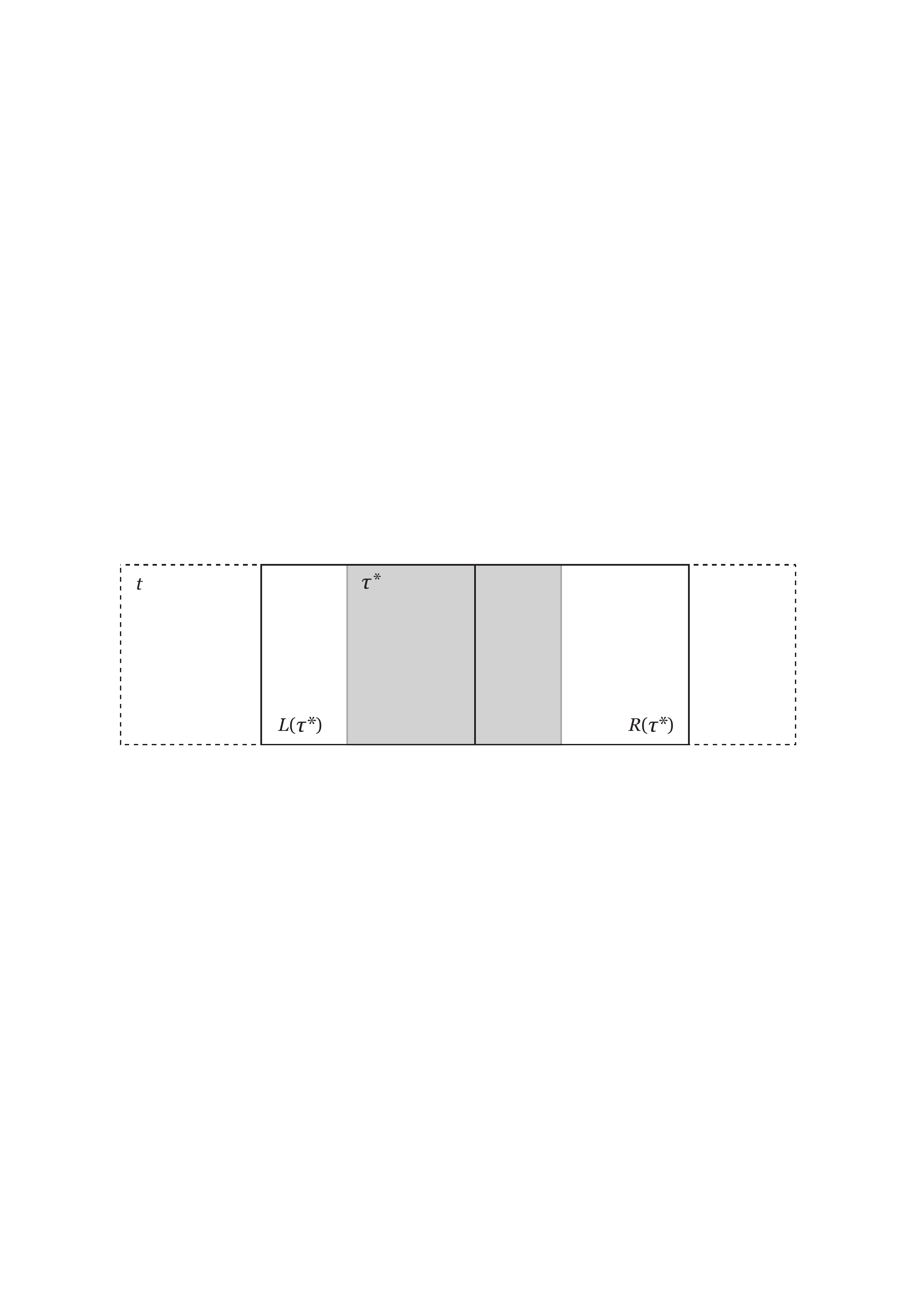}
\end{center}
\caption{The gray rectangle in the middle is $\tau^*$ while the other two are its dyadic translations $L(\tau^*)$ and $R(\tau^*)$. They are inside of $t$, the dashed rectangle.}
\label{fig:dyadicLR}
\end{figure}

Denote the dyadic rectangles by $R(\tau^*)$ and $L(\tau^*)$. 
Note that $|R(\tau^*)|=|L(\tau^*)|=|\tau^*|$ and thus 
\[
\frac{1}{|\tau^*|}\int_{\tau^*} f\, dxdy \leq \frac{1}{|R(\tau^*)|}\int_{R(\tau^*)} f\, dxdy + \frac{1}{|L(\tau^*)|}\int_{L(\tau^*)} f\, dxdy. 
\]
Hence for at least one rectangle $K(\tau^*)\in \{ R(\tau^*),L(\tau^*)\}$ we will have that %(say it is $L(\tau^*)$)
\begin{equation} \label{ineq:L/R}
    \frac{1}{|K(\tau^*)|}\int_{K(\tau^*)} f\, dxdy\geq \frac{n}{2}
\end{equation}
by Lemma \ref{lem:concentration_support}. 
Now, for each $\tau^*$ we consider the union of all such rectangles $K(\tau^*)$ for which \eqref{ineq:L/R} holds. 
Specifically, let 
\[
E'=\bigcup_{F\in \mathcal{A}} \bigcup_{k=0}^{n-1} \bigcup_{t\in \mathcal{T}_k}K(\tau^*).
\]
(Recall here that $\tau^*$ is an abbreviation of $t^*(\omega_t)$ for $t\in \mathcal{T}_k$.)
According to the definition, with the aid of Lemma \ref{lem:pairing}-2), we have that both $R(\tau^*)$ or $L(\tau^*)$ are inside of $t$. 
Hence, due to Lemma \ref{Lmm-cov} we have that $|E'|\geq 1/2$. 
%Set \[N'=\bigcup_{F\in \mathcal{A}} N,\]where $N\subset [0,1)^2$ is defined in Lemma \ref{lem:def_N}. 
Define  
\begin{equation} \label{def:E0}
    E=E'\setminus N, 
\end{equation} 
where $N\subset [0,1)^2$ is defined in Lemma \ref{lem:def_N}. 
Since $|N|<\varepsilon$ by Lemma \ref{lem:def_N}-2), one has $|E|\geq 1/3$. 
We also automatically have the first assertion for every point in $E$. 

Next, we will prove the second statement. 
Since $E\cap N=\emptyset$ and $N$ is defined through the cases studies 1-ii) and 2-ii) in the proof of Lemma \ref{lem:removalfamily}, it is enough to show the second assertion for $R\in \mathcal{F}_{C}$ belonging to the cases 1-i) and 2-i) there. 
Note here that the estimate \eqref{ineq:card} is verified for every $R\in \mathcal{F}_{C}$ with sides considered in the cases 1-i) and 2-i). 
Hence one can obtain the desired estimate 
\[
\frac{1}{|R|} \int_{R}f\,dxdy \leq 3
\]
for every $z\in [0,1)^2\setminus N$ and every $R\in \mathcal{F}_{C}$ with $R\ni z$. 
Proposition \ref{prop:main1} is obtained. 

\begin{remark}\label{rem:RM1}
We remark that instead of the unit square we could do the same constructions inside any dyadic square. For our purposes in Proposition \ref{prop:main2} below it will be more convenient to consider a partition of the unit square into smaller, dyadic squares and carry out the same constructions inside each tiny square. Then, inside each square we will find the corresponding sets $E,N$ satisfying the bounds $|N|\leq |Q|\varepsilon$ and $|E|\geq (1/3)|Q|$. The estimates \eqref{cond:conv0} and \eqref{cond:div0} will hold for rectangles that are strictly inside the partition squares. So to prove the analog of Proposition \ref{prop:main1} in this case, it will remain to take care of rectangles from $\mathcal{F}_C$ that are not entirely contained inside a partition square. For this we can write $R=\cup_{i=1}^N R_i$, where $R_i$ entirely belongs to a partition square. Note that since all the rectangles are dyadic then the rectangles $R_i$ will have identical size. Then the property \eqref{cond:conv0} can be achieved as follows
\[
\frac{1}{|R|}\int_R f\, dxdy=\frac{1}{N}\sum_i\frac{1}{|R_i|}\int_{R_i} f\, dxdy\leq C_0.
\]
\end{remark}

\subsection{}
We now prove an extension of Proposition \ref{prop:main1}.

\begin{proposition} \label{prop:main2}
There is a constant $C_0>0$, so that for every $\varepsilon >0$ and every $M>0$, there exists a function $f\in L^\infty([0,1)^2)$, with $f \geq 0$ and
\begin{equation}\label{norm}
\|f\|_1 \leq 1,
\end{equation}
for which one can find a subset $E\subset [0,1)^2$, with $|E|\geq 1-\varepsilon$, such that for every $z\in E$ there exists a dyadic rectangle $R$ from $\mathcal{F}_D$, such that $z\in R$ and
\begin{equation}\label{abv}
\frac{1}{|R|}\int_{R} f\, dxdy \geq M,
\end{equation}
and for every $z\in E$ and any dyadic rectangle $R \in \mathcal{F}_{C}$, with $z\in R$, we have that 
\begin{equation}\label{bll}
\frac{1}{|R|}\int_{R} f\, dxdy \leq C_0.
\end{equation}
\end{proposition}
\begin{proof}

For every $k\in \mathbb{N}$, let $\varepsilon_k=\varepsilon/2^{k+1}$. 
Let $n_1<n_2<\cdots <n_k<\cdots $ be an increasing sequence of positive integers.
For each $k\in \mathbb{N}$, consider a partition $\mathcal{P}_k$ of the unit square into dyadic squares $Q_{k,j}$ of size $1/2^{n_k}$. 
Inside each square $Q_{k,j}$, we repeat the same procedure as in Proposition \ref{prop:main1} as described in Remark \ref{rem:RM1}. 
Hence we will find positive functions $f_{k,j}\in L^\infty$ with $\| f_{k,j}\|_{L^1} \leq (1+o(1))|Q_{k,j}|$, and $E_{k,j},N_{k,j}\subset Q_{k,j}$ with $|E_{k,j}|\geq (1/3)|Q_{k,j}|$ and $|N_{k,j}|\leq \varepsilon_k |Q_{k,j}|$ such that \eqref{cond:div0} and \eqref{cond:conv0} holds with $n\geq 2\max\{k^3,M^3\}$ which is a power of two. 
Namely, for every $z\in E_{k,j}$ there exists a dyadic rectangle $R\in \mathcal{F}_D$ with $R\ni z$ such that 
\[
\frac{1}{|R|}\int_R f_{k,j}\, dxdy \geq \frac{n}{2}\geq \max \{ k^3,M^3\},
\]
and for every $z\in Q_{k,j}\setminus N_{k,j}$ and any dyadic rectangle $R\in \{R\in \mathcal{F}_{C}\colon R\subset Q_{k,j}\}$ with $R\ni z$, one has
\begin{equation} \label{ineq:ub_Q}
    \frac{1}{|R|}\int_R f_{k,j}\, dxdy \leq 3.
\end{equation}
Define 
\begin{align*}
    E_k&=\bigsqcup_{j=1}^{2^{2n_k}}E_{k,j}, \\
    N_k&=\bigsqcup_{j=1}^{2^{2n_k}}N_{k,j}, 
\end{align*}
and 
\[
f_k=\sum_{j=1}^{2^{2n_k}}f_{k,j}.
\]
Then for each $k\in \mathbb{N}$ one has $|E_k|\geq 1/3$, $|N_k|\leq \varepsilon_k$ and $\|f_k\|_{L^1}\leq 1+o(1)$. 

Now, given $L\in \mathbb{N}$ to be determined later, define
\[
E=\left(\bigcup_{k=1}^{L} E_k\right)\setminus \left(\bigcup_{k=1}^{L} N_k\right),
\]
and 
\[
f=\sum_{k=1}^L \frac{1}{k^2}f_k. 
\]
Then one has
\[
\|f\|_{L^1}\leq \sum_{k=1}^L \frac{\|f_k\|_{L^1}}{k^2}\leq (1+o(1))\sum_{k=1}^L \frac{1}{k^2}\leq C 
\]
for some $C>0$, and for every $z\in E$ there exists a dyadic rectangle $R\in \mathcal{F}_D$ with $R\ni z$ such that
\begin{align*}
    \frac{1}{|R|}\int_R f\, dxdy 
    &\geq \frac{1}{k^2}\frac{1}{|R|}\int_R f_k\, dxdy \\
    &\geq \frac{1}{k^2} \frac{1}{|R|}\int_R f_{k,j}\, dxdy \geq \frac{1}{k^2} \frac{n}{2}\geq \frac{\max \{ k^3,M^3\}}{k^2}\geq M. 
\end{align*}
We also have for every $z\in E$ and any dyadic rectangle $R \in \mathcal{F}_C$, with $z\in R$, 
\[
\frac{1}{|R|}\int_R f\, dxdy\leq 3 \sum_{k=1}^L \frac{1}{k^2}
\]
by \eqref{ineq:ub_Q} with Remark \ref{rem:RM1}. 

It remains to show $|E|>1-\varepsilon$. 
Note first that 
\[
\left| \bigcup_{k=1}^{L} N_k\right| \leq \sum_{k=1}^{L}|N_k|\leq \sum_{k=1}^{L}\varepsilon_k=\sum_{k=1}^{L}\frac{\varepsilon}{2^{k+1}} < \frac{\varepsilon}{2}.
\]
Thus it is enough to show $|\cup_{k=1}^{L} E_k|\geq 1-(\varepsilon/2)$, and this will be achieved by making $n_k$ grow sufficiently fast and taking $L$ sufficiently large. 
For $k\in \mathbb{N}$, we let $E(k)=E_1\cup \dots \cup E_k$. 
Note that if $n_k$ grows fast then the partition at step $k+1$ can be made so small that $E_{k+1}\setminus E(k)$ will fill up almost $1/3$rd of the compliment of $E(k)$. 
Hence, by taking $L$ large enough we can achieve the bound $|E(L)|\geq 1-(\varepsilon/2)$. 
It then follows that 
\[
|E|\geq |E(L)|-\left| \bigcup_{k=1}^{L} N_k\right| > 1-\frac{\varepsilon}{2}-\frac{\varepsilon}{2}=1-\varepsilon.
\]
Proposition \ref{prop:main2} is obtained. 
\end{proof}

\section{Proof of Theorem \ref{thm:main2}}

\begin{proof}
Assume \eqref{cond:dst}. 
We use Proposition \ref{prop:main2} for $\varepsilon_k=1/2^k$ and  $M=k^3$. 
We will get a sequence of functions $f_k\in L^\infty$ and sets $E_k\subset [0,1)^2$ such that $|E_k|>1-\varepsilon_k$. 
We then consider the function
\[
f=\sum_{n=1}^{\infty} \frac{1}{n^{2}} f_{n}.
\]
Note that by Proposition \ref{prop:main2}, \eqref{norm}, we have
\[
\left\|\sum_{n=1}^{\infty} \frac{1}{n^{2}} f_{n}\right\|_{L^{1}} \leq \sum_{n=1}^{\infty} \frac{1}{n^{2}}\left\|f_{n}\right\|_{L^{1}}<\infty .
\]
Hence, $f$ is well defined. It is positive and $f \in L^{1}$. 
By the Borel-Cantelli lemma, we also have that
\[
\textrm{Leb} \left(\liminf _{n \rightarrow \infty} E_{n}\right)=1
\]
Since $f_{k} \in L^{\infty}\subset L^1 \ln^+ L^1$, then by the Jessen-Marcinkiewicz-Zygmund theorem \cite{Jessen-Marcinkiewicz-Zygmund1935}, we have that for all $k \in \mathbb{N}$, there is $\Gamma_{k} \subset [0,1)^{2}$ with $\textrm{Leb}\left(\Gamma_{k}\right)=1$ such that
\begin{equation} \label{cond:convergence}
    \lim _{\substack{\text{diam} R \rightarrow 0; \\z\in R \in \mathcal{R}}} \frac{1}{|R|} \int_{R} f_{k}\, d xdy=f_{k}(z)
\end{equation}
for every $z \in \Gamma_{k}$. 
It follows from \eqref{cond:convergence} and Proposition \ref{prop:main2}-\eqref{bll} that for every $z \in E_{k} \cap \Gamma_{k}$ we have
$$
\left|f_{k}(z)\right| \leq C_0. 
$$
Define $\Gamma_{\infty}=\bigcap_{k=1}^{\infty} \Gamma_{k}$ and
$$
\Lambda=\left(\liminf _{n \rightarrow \infty} E_{n}\right) \cap \Gamma_{\infty}
$$
Clearly $\textrm{Leb}(\Lambda)=1$. 

First, by assumption we have that almost every $z\in [0,1)^{2}$ eventually belongs to all sets $E_k$, i.e.  there exists $K=K(z)\in \mathbb{N}$ so that for all $s\geq K$ we have $z\in E_s\cap \Gamma_{\infty}$. If $z \in E_s$, for some $s$, then by \eqref{abv} we can find $R\in \mathcal{F}_D$ with $R\ni z$ so that 
\[
\frac{1}{|R|}\int_R f\, dxdy \geq \frac{1}{s^2}
\frac{1}{|R|}\int_R f_s\, dxdy\geq \frac{s^3}{s^2}=s.
\]

Next, for $R\in \mathcal{F}_C$ write
\[
\frac{1}{|R|}\int_R\sum_{n=1}^{\infty} \frac{1}{n^{2}} f_n\, dxdy = \sum_{n=1}^N \frac{1}{n^2|R|}\int_R f_n\, dxdy + \sum_{n=N+1}^\infty \frac{1}{n^2|R|}\int_R f_n\, dxdy.
\]
Since almost every $z$ eventually belongs to all sets $E_k$, then for large enough $N$ and property \eqref{bll}
\[
\sum_{n=N+1}^\infty \frac{1}{n^2|R|}\int_R f_n\, dxdy
\leq C_0\sum_{n=N+1}^\infty \frac{1}{n^2}.
\]
This can be made small if $N$ is large. 
While for the  first term we have {by \eqref{cond:convergence}}
\[
\lim_{\substack{\text{diam} R \rightarrow 0; \\z\in R \in \mathcal{F}_C}}\sum_{n=1}^N \frac{1}{n^2|R|}\int_R f_n\, dxdy=\sum_{n=1}^{N} \frac{1}{n^{2}} f_{n}(z)
\]
for $z\in \Lambda$.
Thus
\[
\lim_{\substack{\text{diam} R \rightarrow 0; \\z\in R \in \mathcal{F}_C}}\frac{1}{|R|}\int_R\sum_{n=1}^{\infty} \frac{1}{n^{2}} f_n\, dxdy=\sum_{n=1}^{\infty} \frac{1}{n^{2}} f_{n}(z)=f(z),
\]
for almost every $z$. 

The proof of the opposite direction is analogous to the necessity part of Theorem \ref{main1} below, so we will skip it.
\end{proof}

\section{Proof of Theorem \ref{main1}}

\begin{proof}[Proof of Theorem \ref{main1}]
We first prove sufficiency of \eqref{cond:dst}. 
By Proposition \ref{haar} and the discussion right after it we have for every $z\in [0,1)^2$ that 
\begin{equation} \label{cond:haar_av}
    S_{n,m}f(z)=\frac{1}{|I_{n,m}|}\int_{I_{n,m}}f\, dxdy,
\end{equation}
where $I_{n,m}(z)$ is a dyadic rectangle containing $z$ with sides $1/2^k$ and $1/2^{s}$, or $1/2^{k+1}$ and $1/2^{s+1}$ (see Proposition \ref{haar}). 
(Recall also definitions \eqref{def:C} and \eqref{def:D}). 

Consider the sets $\mathcal{N}$ and $\mathcal{M}$ and the associated bases $\mathcal{F}_{C}$ and $\mathcal{F}_{D}$. If now $n,m \in \mathcal{N}$ or $n,m \in \mathcal{M}$, then the rectangle $I_{n,m}$ above will belong to either $\mathcal{F}_{C}$ or $\mathcal{F}_{C+1}$, or $\mathcal{F}_{D}$ and $\mathcal{F}_{D+1}$, where $C+1=\{c+1\colon c \in C\}$.
Define new sets $\mathbf{C}=C\cup (C+1)$ and $\mathbf{D}=D+1$. 
Note that $\mathbf{C}$ and $\mathbf{D}$ satisfy the assumption \eqref{cond:dst} since it is also satisfied by $C$ and $D$.  
Hence, by Theorem \ref{thm:main2}, there exists a non-negative function $f\in L^1([0,1)^2)$ such that for almost every $z$ we have 
\[
\lim_{\substack{\diam R \rightarrow 0\colon \\ z\in R\in \mathcal{F}_{\mathbf{C}}}}\frac{1}{|R|}\int_{R}f\,dxdy=f(z)
\]
and
\[
\limsup_{\substack{\diam R \rightarrow 0; \\ z\in R\in \mathcal{F}_{\mathbf{D}}}}\frac{1}{|R|}\int_{R}f\, dxdy=\infty.
\]
In view of \eqref{cond:haar_av}, it follows from the first relation that for almost every $z$ we have
\[
\lim_{\substack{n,m\rightarrow \infty; \\ n,m \in C}}S_{n,m}f(z)=\lim_{\substack{n,m\rightarrow \infty; \\n,m \in \mathbf{C}}}S_{n,m}f(z)=f(z).
\]

To see the second part, note that for every $z$ the dyadic rectangle with sides $1/2^{k+1},1/2^{s+1}$ which contains $z$, is also contained in the rectangle with sides $1/2^k,1/2^{s}$ containing $z$. Hence, since $f$ is positive and we have divergence with respect to the dyadic rectangle with sides $1/2^{k+1},1/2^{s+1}$ ($k,s \in C$) then we also have it for the dyadic rectangle with sides $1/2^k,1/2^{s}$. This proves the divergence part of the theorem.

\begin{comment}
   
Thus, for almost every $(x,y)$ there exists a sequence from $\mathcal{B}_0$ such that we have the above. However, we need to show that
$$
\lim_{n_k,m_k \rightarrow \infty; I_{n,m}\in \mathcal{F}_{\mathcal{B}_0}}\frac{1}{|I_{n,m}(x,y)|}\int_{|I_{n,m}(x,y)|}fdt=\infty.
$$
In the same way we can show that
$$
\lim_{n,m \rightarrow \infty; I_{n,m}\in \mathcal{F}_{\mathcal{B}_0}}\frac{1}{|I_{n+1,m+1}(x,y)|}\int_{|I_{n+1,m+1}(x,y)|}fdxdy=\infty.
$$
Clearly, from the sequence $S_{n_k,m_k}f(x,y)$ we can find a subsequence which is a subsequence eith of or. Hence, we will have convergence for both.
\end{comment}

Next, we prove the necessity of \eqref{cond:dst} by contraposition. 
Hence, assume that \eqref{cond:dst} fails. 
Then there exists an integer $d>0$ so that for every $s \in D \cup (D+1)$ we have that
\begin{equation} \label{cond:cntraposition}
    \Big(C\cup (C+1)\Big)\cap [s-d, s+d]\neq \emptyset.
\end{equation}
Suppose 
$$
\lim_{\substack{n,m \rightarrow \infty; \\n,m\in \mathcal{N}}}S_{n,m}f(z)=f(z)
$$
for almost every $z$. 
For a given $\delta>0$ consider the set
\[
E_{\delta}=\left\{z \in \mathbb{R}^{2}\colon \sup _{n,m \geq \frac{1}{\delta}, n,m \in \mathcal{N}}\left|S_{n,m}f(z)-f(z)\right|<1\right\} .
\]

We have $\left|E_{\delta}\right|>0$ for $\delta>0$ small enough. 
Then for the characteristic function $\mathbbm{1}_{E_{\delta}}$ of $E_{\delta}$, we will have by the Jessen-Marcinkiewicz-Zygmund theorem \cite{Jessen-Marcinkiewicz-Zygmund1935}, that for almost all points $z \in E_{\delta}$ are Lebesgue points, namely for almost every $z \in E_{\delta}$ we have
\[
\lim _{\substack{\text{diam} R \rightarrow 0; \\ z \in R \in \mathcal{R}}} \frac{\left|R \cap E_{\delta}\right|}{|R|}=1 .
\]
Let $z \in E_\delta$. 
Assume $n,m\in \mathcal{M}$ are so large that for the interval $B=I_{n,m}(z)$ from \eqref{Haar:rect}, i.e. $S_{n,m}f(z)=(1/|I_{n,m}(z)|)\int_{I_{n,m}(z)}fdxdy$, we have that
\[
\frac{\left|B \cap E_{\delta}\right|}{|B|}>1-c,
\]
where $c\in (0,1)$ is to be chosen later. 
Assume the sides of $B$ are $1/2^s$ and $1/2^k$, where $s,k\in D \cup (D+1)$. 
We now represent $B$ as a union of dyadic rectangles with sides $1/2^{s+d}$ and $1/2^{k+d}$, i.e. $B=\sqcup _{q}B_q$. 
Note that, since $d$ is fixed and the constant $c$ above can be taken arbitrarily small, then for an appropriate choice of $c$ we can make sure that each of the rectangles $B_q$ has a non-empty intersection with $E_\delta$. 
Thus for each $B_q$ we can chose a point $z\in E_\delta \cap B_q$. 
Then by \eqref{cond:cntraposition} and Proposition \ref{haar}, there exists $q=(n_1,m_1)$, with $n_1,m_1\in \mathcal{N}$, a rectangle $I_q=I_q(z)$ containing $z$, with sides $1/2^{\mu (z)}$ and $1/2^{\nu (z)}$, where $\mu (z),\nu (z)\in C \cup (C +1)$ so that $\mu (z)\in [s-d,s+d]$ and $\nu (z)\in [k-d,k+d]$, respectively and
\[
S_{n_1,m_1}f(z)=\frac{1}{|I_q|}\int_{I_q}f\, dxdy<1+f(z).
\]
Thus, we will also have that $B_q \subset I_q$.
Note that %$B_k \subset I_k$ and that 
\[
\frac{|B_q|}{|I_q|}\geq \frac{1/2^{s+d}\cdot 1/2^{k+d}}{1/2^{s-d}\cdot 1/2^{k-d}}=\frac{1}{2^{4d}} .
\] 
%since $B_k$ is dyadic and its sides are less than the sides of $I_k$. 
Repeating the same argument for all remaining dyadic rectangles $B_q$ we can find a collection of rectangles $\{I_q{ \in \mathcal{F}_{C\cup (C+1)}}\}_q$ such that $B \subset \cup_{q} I_q$. 
%Next, by assumption 
Then for almost every $z \in E_\delta$ we have
\begin{align*}
    \int_{B} f\, dxdy 
    &\leq \sum_{q} \int_{I_{q}} f\, dxdy \\ 
    &\leq \left( 1+f(z)\right) \sum_{q}\left|I_{q}\right| \leq \left( 1+f(z)\right) \sum_{q} 2^{4d}|B_q|=\left( 1+f(z)\right) 2^{4d}|B|,
\end{align*}
and hence
\[
\frac{1}{|B|}\int_{B} f\, dxdy \leq \left(1+f(z)\right)2^{4d}.
\]
Which implies that for almost every $z \in E_\delta$ we have
\[
\lim_{\substack{n,m \rightarrow \infty;\\n,m\in \mathcal{M}}}S_{n,m}f(z)\leq\limsup_{\substack{\diam B \rightarrow 0; \\z \in B \in \mathcal{F}_{D\cup (D+1)}}}\frac{1}{|B|}\int_{B} f\, dxdy<\infty.
\]
This contradicts the assumption that the above limsup is unbounded on the bases $\mathcal{F}_D$. 
This finishes the proof.
\begin{comment}
   
\begin{figure}[H]
\begin{center}
\includegraphics[pagebox=artbox,width=6cm]{covering.png}
\end{center}
\caption{$B$ is covered by rectangles of sides $x,y$. The average has to be large at least for one of these rectangles contradicting the assumption.}
\label{fig:Bk}
\end{figure}
\end{comment}
\end{proof}

\appendix    
\section{The Haar wavelet and its properties}\label{app}

For simplicity, we will formulate the multivariate Haar system only in dimension $2$. 
More general formulations can be found in \cite{Alexits1961, Oniani2012}. 
Our presentation follows the notations of \cite{Oniani2012}. 

Let $\mathbb{N}_{0}=\mathbb{N}\cup \{0\}$ be the set of all non-negative integers and $\mathbb{N}^{2}=\mathbb{N}\times \mathbb{N}$.
We recall the definition of the one dimensional Haar system $\left\{h_{m}\right\}_{m \in \mathbb{N}}$ :
$$
h_{1}(x)=1 \quad\left(x \in \mathbb{T}^{1}\right)
$$
and if
$$
m=2^{k}+i \quad\left(k \in \mathbb{N}_{0}, i = 1,\dots, 2^{k}\right)
$$
then
$$
h_{m}(x)= \begin{cases}2^{k / 2} & \text { if } x \in \left(\frac{i-1}{2^{k}}, \frac{2 i-1}{2^{k+1}}\right) \\ -2^{k / 2} & \text { if } x \in \left(\frac{2 i-1}{2^{k+1}}, \frac{i}{2^{k}}\right) \\ 0 & \text { if } x \notin\left[\frac{i-1}{2^{k}}, \frac{i}{2^{k}}\right]\end{cases}
$$
At inner points of discontinuity $h_{m}$ is defined as the mean value of the limits from the right and from the left, and at the endpoints of $\mathbb{T}^{1}$ as the limits from inside of the interval.
The two dimensional Haar system $\left\{h_{(m,n)}\right\}_{(m,n)\in \mathbb{N}^2}$ is defined as follows:
\[
h_{(m,n)}(x,y)=h_{m}(x) \times  h_{n}(y) \quad \left( (x,y) \in \mathbb{T}^{2}\right) .
\]
For $z \in \mathbb{T}^{2}$, let $H(z)$ be the spectrum of the Haar system at $z$, i.e., 
\[
H(z)=\left\{(m,n)\in \mathbb{N}^{2}\colon h_{(m,n)}(z) \neq 0\right\}.
\]
We denote by $\Delta_{k,s}(z)$ the dyadic rectangle with sides $1/2^k$ and $1/2^s$ that contains $z$.

The next property connects rectangular convergence of Fourier-Haar series with the differentiation of integrals with respect to the basis of dyadic rectangles. 
(See, e.g., [\cite{KashSaak1984}, \text { Ch. 3}, \S 1] or [\cite{Alexits1961}, \text { Ch. 1}, \S 6].) 

Below, by $[a,b]$ with $a,b\in \mathbb{N}$ and $a\leq b$ we mean the set $\{a, a+1, \dots , b\}$

\begin{proposition}\label{haar}
Let $f \in L^1\left(\mathbb{T}^{2}\right)$, $z \in \mathbb{T}^{2}$ and $(m,n)\in \mathbb{N}^{2}$. 
Then the following assertions hold: let $m=2^k+i$ and $n=2^s+j$, with $i=1,\dots, 2^k$ and $j=1,\dots ,2^s$;
\begin{enumerate}
    \item[1)] If $H(z) \cap(\left[2^{k}+1, m\right] \times \left[2^{s}+1, n\right])\neq \emptyset$, then
    \[
    S_{m,n}f(z)=S_{2^{k+1},2^{s+1}}f(z)=\frac{1}{|\Delta_{k+1,s+1}(z)|}\int_{\Delta_{k+1,s+1}(z)}f\, dxdy.
    \]
    \item[2)] 
    If $H(z) \cap(\left[2^{k}+1, m\right] \times \left[2^{s}+1, n\right])= \emptyset$, then
    \[
    S_{m,n}f(z)=S_{2^{k},2^{s}}f(z)=\frac{1}{|\Delta_{k,s}(z)|}\int_{\Delta_{k,s}(z)}f\, dxdy.
    \]
\end{enumerate}
\end{proposition}

In other words for each $z\in \mathbb{T}^{2}$ and every $m,n\in \mathbb{N}$, we have that
\begin{equation}\label{Haar:rect}
S_{m,n}f(z)=\frac{1}{|I_{m,n}(z)|}\int_{I_{m,n}(z)}f\ dxdy
\end{equation}
where $I_{m,n}(z)$ is a dyadic rectangle with sides $1/2^k$ and $1/2^s$ or $1/2^{k+1}$ and $1/2^{s+1}$ containing $z$. 

\begin{comment}
   
For $n,m$ let $n=2^k+i$ and $m=2^s+j$, where $k,s\geq 1$ and $i< 2^k$, $j< 2^s$. It is know that

$$
S_{n,m}(x,y)= \begin{cases}S_{2^k,2^s}(x,y) & \text { if } x \in\left[0, \frac{i}{2^{k}}\right) \text { or }y \in\left[0, \frac{j}{2^{s}}\right)\\ S_{2^{k+1},2^{s+1}}(x,y) & \text { if } x \in\left(\frac{i}{2^{k}},1\right) \text { or }y \in\left(\frac{j}{2^{s}},1\right)\end{cases}
$$
It is also known that 
$$
S_{2^{k+1},2^{s+1}}(x,y)= \begin{cases}\frac{1}{|I_{i,+}||I_{j,+}|}\int_{I_{i,+}\times I_{j,+}}f(t_1,t_2)dt_1 dt_2 & \text { if } (x,y)\in I_{i,+}\times I_{j,+}\\
\frac{1}{|I_{i,-}||I_{j,+}|}\int_{I_{i,-}\times I_{j,+}}f(t_1,t_2)dt_1 dt_2 & \text { if } (x,y)\in I_{i,-}\times I_{j,+}\\ 
\frac{1}{|I_{i,+}||I_{j,-}|}\int_{I_{i,+}\times I_{j,-}}f(t_1,t_2)dt_1 dt_2 & \text { if } (x,y)\in I_{i,+}\times I_{j,-}\\ 
\frac{1}{|I_{i,-}||I_{j,-}|}\int_{I_{i,-}\times I_{j,-}}f(t_1,t_2)dt_1 dt_2 & \text { if } (x,y)\in I_{i,-}\times I_{j,-}\end{cases}
$$
\end{comment}

\section*{Acknowledgements}
The authors would like to thank Shigeki Akiyama and Tomas Persson for helpfull discussions. And also to \text { Håkan Hedenmalm } for usefull suggestions. The first author is partially supported by Japan Society for the Promotion of Science (JSPS)
KAKENHI Grant Number 19K03558. The second author is supported by the Knut and Alice
Wallenberg foundation of Sweden (KAW) .

\end{document}